
%
%
\input amstex	             
\documentstyle{amsppt}	     
\NoBlackBoxes	             
\magnification=\magstep1

\topmatter
\title
  On polarized manifolds of sectional genus three
\endtitle
\author
  Hironobu Ishihara
\endauthor
\affil
  Department of Mathematics,
  Faculty of Science,
  Tokyo Institute of Technology
\endaffil
\address
  Oh-okayama,Meguro,Tokyo,152,Japan
\endaddress
\email
  ishihara\@math.titech.ac.jp
\endemail
\endtopmatter

\document

\subhead
Introduction
\endsubhead

Let $L$ be an ample line bundle on a complex projective manifold $M$ of
dimension $n\geq 2$.
The sectional genus $g=g(M,L)$ of a polarized manifold
$(M,L)$ is defined by the formula $2g(M,L)-2=(K+(n-1)L)L^{n-1}$,
where $K$ is the canonical bundle of $M$.
For polarized manifolds over $\Bbb C$,it is known that $g$ takes
non-negative integers (\cite{F6;Corollary 1}).

In many papers the structure of $(M,L)$ with low $g$ has been studied:
see \cite{F6} for $g\leq 1$ ; \cite{BeLP} for $g=n=2$ ; \cite{F7} for $g=2$ ;
\cite{Ma} for $g=3$ and $n=2$.
In this paper we study the case $g=3$ and $n\geq 3$.
Under the additional condition that $L$ is spanned,the classification was
partially known by \cite{BiLL}.
Here we study $(M,L)$ without this hypothesis.

This paper is organized as follows.
In \S1 we show that $(M,L)$ with $g=3$ and $n\geq 3$ is
one of the following types.

\roster
\item  There is an effective divisor $E$ on $M$ such that
       $(E,L_E)\simeq (\Bbb P^{n-1},\Cal O(1))$ and
       $[E]_E=\Cal O(-1)$.

\item  There is a fibration $\varPhi:M\to C$ over a smooth curve $C$ such that
       $(F,L_F)\simeq (\Bbb P^2,\Cal O(2))$ for every fiber $F$ of $\varPhi $.

\item  There is a fibration $\varPhi:M\to C$ over a smooth curve $C$ such that
       every fiber $F$ of $\varPhi$ is a hyperquadric in $\Bbb P^n$
       and $L_F=\Cal O(1)$.

\item  $(M,L)$ is a scroll over a smooth surface.

\item  $K+(n-2)L$ is nef.

\item  $(M,L)$ is a scroll over a smooth curve of genus three.
\endroster

\flushpar
In \S2 we study the case \thetag{4},in \S3 we study the case \thetag{3},
in \S4 we study the case \thetag{2},
and in \S5 we study the cases \thetag{1} and \thetag{5}.
Although our results are far from being complete,
they are very similar to those in case $g=2$.

The author would like to express his sincere thanks to Professor T.~Fujita
for kind encouragement and for many valuable comments
during the preparation of this paper.

\subhead
Notation
\endsubhead

Basically we use the customary notation in algebraic geometry as in \cite{H2}.
All varieties are defined over $\Bbb C$ and assumed to be complete.
Vector bundles are often identified with locally free sheaves
of their sections,and these words are used interchangeably.
Line bundles are identified with linear equivalence classes
of Cartier divisors,and their tensor products are denoted additively,
while we use multiplicative notation for intersection products in Chow rings.
The numerical equivalence of line bundles is denoted $\equiv$,
while we use $=$ for linear equivalence.
The linear equivalence class is denoted by $[~]$,and its corresponding
invertible sheaf is denoted by $\Cal O[~]$.

Given a morphism $f:X\to Y$ and a line bundle $A$ on $Y$,
we denote $f^*A$ by $A_X$,or sometimes by $A$ for short
when there is no danger of confusion.
The canonical bundle of a manifold $M$ is denoted by $K^M$,
unlike the customary notation $K_M$.
The $\Cal O(1)$'s of projective spaces
$\Bbb P_\alpha,\Bbb P_\beta,\dots$ will be denoted by
$H_\alpha,H_\beta,\dots$.
Given a vector bundle $\Cal E$ on $X$,we denote by
$\Bbb P_X(\Cal E)$ (or $\Bbb P(\Cal E)$)
the associated projective space bundle,and denote by $H(\Cal E)$
the tautological line bundle on $\Bbb P(\Cal E)$ in the sense of
\cite{H2}.
The pair $(\Bbb P(\Cal E),H(\Cal E))$ is called the scroll of $\Cal E$.

\subhead
\S1  Classification;~first step
\endsubhead

Throughout this paper a polarized manifold $(M,L)$ is a pair of
a nonsingular projective variety $M$ over $\Bbb C$ and an ample line
bundle $L$ on $M$.
We consider the case with $n=\dim M\geq 3$ and denote by $K$
the canonical bundle of $M$.
First we review known results about polarized manifolds.

\proclaim{{\rm(1.1)}Theorem({\rm \cite{F6;Theorem 1}})}
Let $(M,L)$ be a polarized manifold.
Then $K+nL$ is nef unless $(M,L)\simeq(\Bbb P^n,\Cal O(1))$.
\endproclaim

\proclaim{{\rm(1.2)}Theorem({\rm \cite{F6;Theorem 2}})}
Let $(M,L)$ be a polarized manifold with $n\geq 2$.
Suppose that $K+nL$ is nef.
Then $K+(n-1)L$ is nef unless either
\roster
\item"({\rm a})" $M$ is a hyperquadric in $\Bbb P^{n+1}$ and $L=\Cal O(1)$,
\item"({\rm b})" $(M,L)\simeq(\Bbb P^2,\Cal O(2))$,or
\item"({\rm c})" $(M,L)$ is a scroll over a smooth curve of genus $g(M,L)$.
\endroster
\endproclaim

We denote by $g=g(M,L)$ the sectional genus of $(M,L)$.
In the above cases (a) and (b),we have $g=0$.
Thus we obtain

\proclaim{{\rm(1.3)}Corollary}
For a polarized manifold $(M,L)$ with $g=3$ and $n\geq 2$,
if $K+(n-1)L$ is not nef,then $(M,L)$ is a scroll over a smooth curve
of genus three.
\endproclaim

When $K+(n-1)L$ is nef,we use the following theorem.

\proclaim{{\rm(1.4)}Theorem({\rm \cite{F6;Theorem 3}})}
Let $(M,L)$ be a polarized manifold with $n\geq 3$.
Suppose that $K+(n-1)L$ is nef.
Then $K+(n-2)L$ is nef except the following cases.
\roster
\widestnumber\item{(b1-Q)}
\item"({\rm a})" There is an effective divisor $E$ on $M$ such that
                 $(E,L_E)\simeq(\Bbb P^{n-1},\Cal O(1))$ and
                 $[E]_E=\Cal O(-1)$.
\item"({\rm b0})" $(M,L)$ is a Del Pezzo manifold~(i.e.~$K+(n-1)L=0$),
                  $(\Bbb P^3,\Cal O(2)),(\Bbb P^3,\Cal O(3))$,
                  \linebreak $(\Bbb P^4,\Cal O(2))$
                  ,or a hyperquadric in $\Bbb P^4$ with $L=\Cal O(2)$.
\item"({\rm b1})" There is a fibration $\varPhi:M\to C$ over a smooth
                  curve $C$ with one of the following properties:
\item"({\rm b1-V})" $(F,L_F)\simeq(\Bbb P^2,\Cal O(2))$ for every fiber $F$ of
                    $\varPhi$;
\item"({\rm b1-Q})" every fiber $F$ of $\varPhi$ is a hyperquadric
                    in $\Bbb P^n$ and $L_F=\Cal O(1)$.
\item"({\rm b2})" $(M,L)$ is a scroll over a smooth surface.
\endroster
\endproclaim

In the case of (b0),we have $g\not=3$.
Thus when $g=3$,Theorem (1.4) is rephrased as below.

\proclaim{{\rm(1.5)}Theorem}
For a polarized manifold $(M,L)$ with $g=3$ and $n\geq 3$,
if $K+(n-1)L$ is nef,then $(M,L)$ is one of the following types.
\roster
\item There is an effective divisor $E$ on $M$ such that
      $(E,L_E)\simeq(\Bbb P^{n-1},\Cal O(1))$
      and $[E]_E=\Cal O(-1)$.
\item There is a fibration $\varPhi:M\to C$ over a smooth curve
      such that $(F,L_F)\simeq(\Bbb P^2,\Cal O(2))$
      for every fiber $F$ of $\varPhi$.
\item There is a fibration $\varPhi:M\to C$ over a smooth curve
      such that every fiber $F$ of $\varPhi$ is a hyperquadric
      in $\Bbb P^n$ and $L_F=\Cal O(1)$.
\item $(M,L)$ is a scroll over a smooth surface.
\item $K+(n-2)L$ is nef.
\endroster
\endproclaim

We study the above cases in the following sections;
in \S2 we study the case (4),
in \S3 we study the case (3),
in \S4 we study the case (2),and
in \S5 we study the cases (1) and (5).

\subhead
\S2  The case of a scroll over a surface
\endsubhead

In this section we study the case (4) of the theorem (1.5),
following the idea in \cite{F8;\S2}.
{}From the definition of scrolls,we have
$(M,L)\simeq(\Bbb P_S(\Cal E),H(\Cal E))$
for some ample vector bundle $\Cal E$ on a smooth surface $S$.

(2.1)
Since $\Cal E$ is ample,$A:=\det\Cal E$ is ample and
$(S,A)$ is a polarized surface.
A simple computation shows $g(S,A)=g(M,L)=3$,thus the classification is
reduced to the classification of polarized surfaces with $g=3$.

(2.2)
We first recall the definition of the minimalization of polarized surfaces
(For details we refer to \cite{F;\S 14}.).
Let $(S,A)$ be a polarized surface.
For a $(-1)$-curve $E$ on $S$,let $\pi:S\to S^-$ be the contraction of $E$.
Then $A+mE=\pi^*A^-$ for an ample line bundle $A^-$ on $S^-$ and
$m:=AE$ is called the weight of the contraction $\pi:(S,A)\to(S^-,A^-)$.
$\pi$ is said to be admissible if $AZ\geq m$ for any $(-1)$-curve $Z$
on $S$.
After a finite sequence of admissible contractions:
$(S,A)=(S_0,A_0)\overset\pi_1\to\longrightarrow
       (S_1,A_1)\overset\pi_2\to\longrightarrow\dots
                \overset\pi_r\to\longrightarrow
       (S_r,A_r)=(S',A')$,
we obtain that either $(S',A')$ is a $\Bbb P^1$-bundle over a curve
or the canonical bundle $K'$ of $S'$ is nef.
$(S',A')$ is called an admissible minimalization of $(S,A)$.
We stop when $S'\simeq\varSigma_1$ although there is a $(-1)$-curve on $S'$.
The weight sequence of this admissible minimalization is defined to be
$\frak m:=(m_r,\dots,m_1)$,
where $m_j(1\leq j\leq r)$ is the weight of $\pi_j$.
$\frak m$ is known to be an invariant of $(S,A)$
and is independent of the choice of the minimalization process.

Polarized surfaces with $g=3$ are classified in \cite{Ma}.

\proclaim{{\rm (2.3)}Theorem({\rm cf. \cite{Ma}})}
Let $(S,A)$ be a polarized surface.
Taking an admissible minimalization of $(S,A)$:
$(S,A)=(S_0,A_0)\to(S_1,A_1)\to\dots\to(S_r,A_r)=(S',A')$,
we denote by $\frak m=(m_r,\dots,m_1)$ its weight sequence.
We put $(S',A')=(S,A)$ when we need not take a minimalization.
Assume that $g(S,A)=3$ and $A=\det\Cal E$ for some ample vector bundle
$\Cal E$ with $\text{rank}~\Cal E\geq 2$.
Then $(S,A)$ is one of the following types.
\roster
\widestnumber\item{Remark.}
\item"({\rm I})" $K\equiv A$ and $A^2=2$.
\item"({\rm II})" $S$ is a minimal surface of general type,
                  $K^2=1,KA=2$,and $A^2=2$.
\item"({\rm III})" $S$ is a minimal elliptic surface and
                   $(A^2,KA)=(2,2)$ or $(3,1)$.
\item"({\rm IV})" $K'\equiv 0$.
                  $A^2,A^{\prime2}$ and $\frak m$ are as follows:
  $$\matrix
              & A^2 & A^{\prime2} & \frak m=(m_r,\dots,m_1) \\
           1) &  4                                          \\
           2) &  2  &  6  &  (2)
  \endmatrix$$
\item"({\rm V})" There is a vector bundle $\Cal F$ of rank two
                  on an elliptic curve $C$ such that
                  $S'\simeq\Bbb P_C(\Cal F)$,
                  $A'=xH(\Cal F)+B_S$ for some $x\in\Bbb Z$ and
                  $B\in\text{Pic}(C)$.
                  We put $e=c_1(\Cal F)$ and $y=\deg B$.
                  Then $A^2,e,x,y$ and $\frak m$ are as follows.
  $$\matrix
              & A^2 &  e  &  x &   y  & \frak m=(m_r,\dots,m_1) \\
           1) &  8  & 0,1 &  2 &  2-e                           \\
           2) &  6  &  0  &  3 &   1                            \\
           3) &  5  &  1  &  5 &  -2                            \\
           4) &  4  & 0,1 &  4 & 1-2e & (2)                     \\
           5) &  3  & 0,1 &  6 & 1-3e & (3)                     \\
           6) &  3  &  1  &  7 &  -3  & (2)                     \\
           7) &  2  &  0  &  9 &   1  & (4)                     \\
           8) &  2  &  1  & 11 &  -5  & (3)                     \\
           9) &  2  &  0  &  5 &   1  & (2,2)
  \endmatrix$$
\item"Remark." Although $\Cal F$ is normalized in \cite{Ma},
              here we choose $\Cal F$ satisfying
              $c_1(\Cal F)=0$ or $1$
              by tensoring some line bundle.
\item"({\rm VI})" $(S,A)\simeq(\Bbb P^2,\Cal O(4))$.
\item"({\rm VII})" $S'\simeq\varSigma_e:=\Bbb P(\Cal O_{\Bbb P^1}
                                         \oplus \Cal O_{\Bbb P^1}(e))$.
                   We denote by $H$ the tautological bundle on
                   $\varSigma_e$ and by $H_\xi$ the pullback of
                   $\Cal O_{\Bbb P^1}$ on $\varSigma_e$.Then we have
                   $A'=xH+yH_\xi$,
                   where $x$ and $y$ are integers and
                   $E$ is the exceptional curve on $\varSigma_e$.
                   $A^2,x,y,r$,and $\frak m$ are as follows:
  $$\matrix
              & A^2 &   e   & x &  y   &  r & \frak m=(m_r,\dots,m_1) \\
           1) & 16  & 0,1,2 & 2 & 4-e  &  0                           \\
           2) &  4  &  0,1  & 4 & 5-2e &  9 & (2,\dots,2)             \\
           3) &  3  &  0,1  & 6 & 7-2e &  9 & (3,\dots,3)
  \endmatrix$$
\item"({\rm VIII})" There is an integer $j(0\leq j\leq r)$
                  such that $(S_j,-K_j)$ is a Del Pezzo surface,
                  where $K_j$ is the canonical bundle of $S_j$,
                  and $A_j=-aK_j$ for some integer $a$.
                  $A^2,K_j^2,a$,and $(m_j,\dots,m_1)$
                  are as follows:
  $$\matrix
              & A^2 & K_j^2 & a & (m_j,\dots,m_1) \\
           1) &  8  &   2   & 2                   \\
           2) &  5  &   1   & 3 & (2)             \\
           3) &  3  &   1   & 4 & (3,2)           \\
           4) &  2  &   2   & 5 & (4,4,4)         \\
           5) &  2  &   2   & 3 & (2,2,2,2)       \\
           6) &  2  &   1   & 6 & (5,3)
  \endmatrix$$
\endroster
\endproclaim

\demo{Proof} The assertion easily follows from \cite{Ma} and
             the next lemma. \qed
\enddemo

\proclaim{{\rm(2.4)}Lemma}
Let $(S,A)$ be a polarized surface and $\Cal E$ a vector bundle
of rank $n-1$ on $S$ such that
$(M,L)\simeq(\Bbb P_S(\Cal E),H(\Cal E))$ and $\det\Cal E=A$.
Then $A^2=L^n+c_2(\Cal E)\geq 2$ and
$AZ\geq\text{rank}~\Cal E\geq 2$ for any rational curve $Z$ on $S$.
\endproclaim

For a proof of this lemma,see \cite{F8;(2.2)\&(1.3)}.

{}From now on,for some types of $(S,A)$ in the above list,
we would like to classify ample vector bundles $\Cal E$
such that $\det\Cal E=A$.

(2.5)
Suppose that $(S,A)$ is of the type (2.3;I).
{}From (2.4) we have $L^n=c_2(\Cal E)=1$.
Since $K\equiv A$ is ample,$S$ is a minimal surface of general type.
By \cite{Bo;Theorem 9},$p_g=h^0(S,K)\leq 3$ and we have
the following possibilities:
\roster
\item"a)" when $p_g=0,q=h^1(S,\Cal O_S)=0$ by e.g.
               \cite{Be;Th\'eor\`eme X.4};
\item"b)" when $p_g=1,q\leq 1$ by e.g. \cite{Be;Th\'eor\`eme X.4};
\item"c)" when $p_g=2,q=0$ by \cite{Bo;Theorem 12};
\item"d)" when $p_g=3,q=0$ by \cite{Bo;Theorem 10}.
\endroster

(2.6)
Suppose that $(S,A)$ is of the type (2.3;II).
{}From (2.4) we have $L^n=c_2(\Cal E)=1$.
Since $K^2=1$,we have $p_g\leq 2$ by \cite{Bo;Theorem 9} and
                      $q=0$ by \cite{Bo;Lemma 14}.
Moreover,we can rule out the possibility that $p_g=2$.
Assume that $p_g=2$.Then $\chi(\Cal O_S)=3$,where $\chi$ is
the Euler characteristic.
{}From this,we obtain $\chi(A-K)=2$ by the Riemann-Roch theorem.
Since $A(A-K)=0$,we have $h^0(A-K)=0$,thus $h^2(A-K)>0$ from
$\chi(A-K)=2$.
This means $h^0(2K-A)>0$ by Serre duality.
Since $K(2K-A)=0$ and $A(2K-A)=2$,any member of $|2K-A|$
is one $(-2)$-curve.
Hence we have $h^0(2K-A)\leq 1$.
This means $h^2(A-K)\leq 1$,thus $\chi(A-K)\leq 1$.
This is a contradiction,hence $p_g\leq 1$.

(2.7)
Suppose that $(S,A)$ is of the type (2.3;V--1).
We treat this case similarly as in \cite{F8;(2.4)\&(2.5)}.
For every fiber $F$ of the bundle map $\rho:\Bbb P_C(\Cal F)\to C$,
we have $F\simeq\Bbb P^1$ and $AF=2$,hence
rank~$\Cal E=2$ and $\Cal E_F\simeq\Cal O(1)\oplus\Cal O(1)$.
Then $\Cal G:=\rho_*(\Cal E\otimes[-H(\Cal F)])$ is a locally free sheaf
of rank two on $C$ and
$\rho^*\Cal G\simeq\Cal E\otimes[-H(\Cal F)]$.
Moreover we have $M=\Bbb P_S(\Cal E)\simeq S\times_C\Bbb P(\Cal G),
                  c_1(\Cal G)=\deg B=2-e,
                  c_2(\Cal E)=2$,and $L^3=6$.

(2.7.1)
When $e=0$,both $\Cal F$ and $\Cal G$ are semistable.
In fact,for any quotient line bundle $Q$ of $\Cal F$,
denote by $Z$ the section of $\rho$ corresponding to $Q$.
Then $c_1(Q)\geq 0$ since $0<AZ=(2H(\Cal F)+B_S)Z=2c_1(Q)+2$.
Thus $\Cal F$ is semistable.
If $\Cal G$ is not semistable,there exists a quotient line bundle
$Q$ of $\Cal G$ such that $2c_1(Q)<c_1(\Cal G)$.
Then $\rho^*Q\otimes H(\Cal F)$ is a quotient line bundle of $\Cal E$.
Hence $\rho^*Q\otimes H(\Cal F)$ is ample and
$0<c_1^2(\rho^*Q\otimes H(\Cal F))=2c_1(Q)<c_1(\Cal G)=2$.
This is a contradiction,thus $\Cal G$ is semistable.
When $e=1$,the semistability of $\Cal F$ and $\Cal G$ is uncertain.

(2.7.2)
Conversely,let $\Cal F$ and $\Cal G$ be semistable vector bundles
of rank two on $C$ with the property that
$(c_1(\Cal F),c_1(\Cal G))=(0,2)$ or (1,1).
We put $\Cal E=\rho^*\Cal G\otimes H(\Cal F)$,where
$\rho:\Bbb P_C(\Cal F)\to C$ is the bundle map.
Then $\Cal E$ is an ample vector bundle on $S:=\Bbb P_C(\Cal F)$
and a polarized surface $(S,\det E)$ satisfies the condition of (2.3;V--1).

To see this,the ampleness of $\Cal E$ is the only non-trivial part.
Let $F_1$ be any fiber of $\rho:\Bbb P_C(\Cal F)\to C$
and $F_2$ any fiber of $\Bbb P_C(\Cal G)\to C$.
By the semistability criterion in \cite{Mi;(3.1)},
$2H(\Cal F)-eF_1$ and $2H(\Cal G)-(2-e)F_2$ are nef,
where $e=c_1(\Cal F)$.
Since $M:=\Bbb P_S(\Cal E)\simeq S\times_C\Bbb P(\Cal G)$
and $H(\Cal E)=[H(\Cal F)]_M+[H(\Cal G)]_M$,
for the fiber $F:=[F_1]_M=[F_2]_M$ of the morphism $M\to C$,
$2H(\Cal E)-2F=[2H(\Cal F)-eF_1]_M+[2H(\Cal G)-(2-e)F_2]_M$
is nef on $M$.
Hence $H(\Cal E)$ is ample and then $\Cal E$ is ample.

(2.8)
Suppose that $(S,A)$ is of the type (2.3;V--2).
Our argument is similar to (2.7).
For any fiber $F$ of the bundle map $\rho:\Bbb P_C(\Cal F)\to C$,
we have $F\simeq\Bbb P^1$ and $AF=3$,
hence there are only two possibilities:
  $$ \text{a)~rank}~\Cal E=3~\text{and}~\Cal E_F=\Cal O(1)^{\oplus3};
     \text{b)~rank}~\Cal E=2~\text{and}~\Cal E_F=\Cal O(1)\oplus\Cal O(2).
  $$

(2.8.1)
In the case (2.8;a),
$\Cal G:=\rho_*(\Cal E\otimes[-H(\Cal F)])$
is a locally free sheaf of rank three on $C$.
Moreover we have
  $$ \align
    &\Cal E\simeq\rho^*\Cal G\otimes H(\Cal F),
     M=\Bbb P_S(\Cal E)\simeq S\times_C\Bbb P(\Cal F),
     c_1(\Cal G)=1,c_2(\Cal E)=2,L^4=4, \\
    &\Cal F~\text{is semistable and}~\Cal G~\text{is stable.}
     \endalign
  $$
Conversely,let $\Cal F$ and $\Cal G$ be semistable vector bundles on $C$
with the property that rank~$\Cal F=2,c_1(\Cal F)=0$,
                       rank~$\Cal G=3,c_1(\Cal G)=1$.
We put $S=\Bbb P_C(\Cal F)$ and $\Cal E=\Cal G_S\otimes H(\Cal F)$.
Then a polarized surface $(S,\det\Cal E)$ satisfies the condition
of (2.3;V--2).

(2.8.2)
In the case (2.8;b),$\Cal G:=\rho_*(\Cal E\otimes[-2H(\Cal F)])$
is an invertible sheaf on $C$.
Using a natural map $\rho^*\Cal G\to\Cal E\otimes[-2H(\Cal F)]$,
we obtain an exact sequence
$0\to\rho^*\Cal G\otimes[2H(\Cal F)]\to\Cal E\to Q\to 0$
for some line bundle $Q$ on $S$.
Since $\det\Cal E=2H(\Cal F)+\rho^*\Cal G+Q$,
we have $Q=H(\Cal F)+\rho^*T$ for some line bundle $T$ on $C$
with $\deg\Cal G+\deg T=1$.
$Q$ is ample since it is a quotient bundle of $\Cal E$,
hence we have $\deg T>0$.
On the other hand
$c_2(\Cal E)=c_1(\rho^*\Cal G+2H(\Cal F)c_1(H(\Cal F)+\rho^*T)
            =1+\deg T$,
thus we have $\deg T<5$ from (2.4).
Hence there are only the following possibilities:
$$ \alignat 4
 \text{a)}&\deg T=1, &~\deg\Cal G&=0, &~c_2(\Cal E)&=2, &~\text{and}~&L^3=4; \\
 \text{b)}&\deg T=2, &~\deg\Cal G&=-1, &~c_2(\Cal E)&=3, &~\text{and}~&L^3=3;\\
 \text{c)}&\deg T=3, &~\deg\Cal G&=-2, &~c_2(\Cal E)&=4, &~\text{and}~&L^3=2;\\
 \text{d)}&\deg T=4, &~\deg\Cal G&=-3, &~c_2(\Cal E)&=5, &~\text{and}~&L^3=1.
\endalignat $$

(2.9)
Suppose that $(S,A)$ is of the type (2.3;VI).
Our results are similar to \cite{BiLL;(1.4.2)}.
Since $Al=4$ for any line $l$ on $\Bbb P^2$,we have rank~$\Cal E\leq 4$.

(2.9.1)
When rank~$\Cal E=4$,we can prove that
$\Cal E\simeq\Cal O(1)^{\oplus4}$
by similar argument as in \cite{V};see \cite{OSS;Chapter I,(3.2.1)}
for a proof.

(2.9.2)
When rank~$\Cal E=3$,$\Cal E$ is a uniform vector bundle
on $\Bbb P^2$,thus the result \cite{E} applies.
In particular,we have
$\Cal E\simeq\Cal O(1)^{\oplus2}\oplus\Cal O(2)$ or
$T_{\Bbb P^2}\oplus\Cal O(1)$,where $T_{\Bbb P^2}$ is the
tangent bundle of $\Bbb P^2$.

(2.9.3)
When rank~$\Cal E=2$ and $\Cal E$ is a Fano bundle (i.e.
the anti-canonical bundle $-K^{\Bbb P(\Cal E)}$ of
$\Bbb P(\Cal E)$ is ample),the theorem in \cite{SW} applies and
$\Cal E$ is one of the following types.
\roster
\item"a)" $\Cal E\simeq\Cal O(1)\oplus\Cal O(3)$ or
                      $\Cal O(2)\oplus\Cal O(2)$.
\item"b)" There is an exact sequence
          $0\to\Cal O(2)\to\Cal E\to\Cal I_x(2)\to 0$,
          where $\Cal I_x$ is the ideal sheaf of one point
          $x\in\Bbb P^2$.
\item"c)" $\Cal E$ is stable with $c_2(\Cal E)=6$,
          $\Cal E(-1)$ is spanned,and
          $0\to\Cal O(-1)^{\oplus2}\to\Cal O^{\oplus4}
          \to\Cal E(-1)\to 0$ is exact.
\item"d)" $\Cal E$ is stable with $c_2(\Cal E)=7$,
          $\Cal E(-1)$ is spanned,and
          $0\to\Cal O(-2)\to\Cal O^{\oplus3}\to\Cal E(-1)\to 0$
          is exact.
\endroster
Even in the case that $\Cal E$ is not a Fano bundle,
we can apply the argument in \cite{SW}.
As a result,$\Cal E$ is of the type b) above if $\Cal E$ is not stable;
$7\leq c_2(\Cal E)\leq 15$ if $\Cal E$ is stable.

(2.10)
Suppose that $(S,A)$ is of the type (2.3;VII--1).
Since $AF=2$ for every fiber $F$ of $\rho:\varSigma_e\to\Bbb P^1$,
we have rank~$\Cal E=2$ and $\Cal E_F=\Cal O(1)\oplus\Cal O(1)$.
Then $\Cal G:=\rho_*(\Cal E\otimes[-H])$ is
a locally free sheaf of rank two on $\Bbb P^1$ and
$\rho^*\Cal G\simeq\Cal E\otimes[-H]$.
Let $Z$ be the section of $\rho$ corresponding to a natural surjection
$\Cal O_{\Bbb P^1}\oplus\Cal O_{\Bbb P^1}(e)\to\Cal O_{\Bbb P^1}$.
Then $H_Z=\Cal O_Z$ and $[\rho^*\Cal G]_Z\simeq\Cal E_Z$.
Hence $\Cal G$ is ample and we obtain
\roster
\item"a)" $\Cal E\simeq [H+H_{\xi}]\oplus[H+3H_{\xi}]$ or
                       $[H+2H_{\xi}]^{\oplus2}$ when $e=0$;
\item"b)" $\Cal E\simeq [H+H_{\xi}]\oplus[H+2H_{\xi}]$ when $e=1$;
\item"c)" $\Cal E\simeq [H+H_{\xi}]^{\oplus2}$ when $e=2$.
\endroster
In these cases,$c_2(\Cal E)=4$ and $L^3=12$.

\subhead
\S3  The case of a hyperquadric fibration over a curve
\endsubhead

In this section we study the case (3) of the theorem (1.5),
following the idea in \cite{F7;\S3}.

(3.1)
Since $h^0(F,L_F)=n+1$ for every fiber $F$ of $\varPhi$,
$\Cal E:=\varPhi_*\Cal O_M[L]$ is a locally free sheaf
of rank $n+1$ on $C$ and a natural map $\varPhi^*\Cal E\to L$
is surjective.
This yields a $C$-morphism $\rho:M\to\Bbb P_C(\Cal E)$
and for every point $x$ on $C$ the restriction of $\rho$ to
$F_x:=\varPhi^{-1}(x)$ is an embedding of $F_x$ into $\Bbb P^n$.
Hence $\rho$ itself is an embedding and $M$ is a member of
$|2H(\Cal E)+B_{\Bbb P(\Cal E)}|$ for some line bundle $B$ on $C$.
We put $d=L^n,e=c_1(\Cal E),b=\deg B$
and denote by $g(C)$ the genus of $C$.
After simple computation,we get $d=2e+b,2g(C)+e+b=4$,and
$s:=2e+(n+1)b\geq 0$.
Furthermore in the last inequality,equality holds if and only if
every fiber of $\rho$ is smooth.
{}From these results,we have $(n+1)d+s+4ng(C)=8n$,
hence $g(C)=0$ or 1.

(3.2)
We first study the case $g(C)=1$.
In this case,$C$ is an elliptic curve and we have
$e=d-2$ and $b=4-d$ from the equality above.
Hence we obtain $d\leq 6$ since $s\geq 0$ and $n\geq 3$.

(3.3)
We study the ampleness of $\Cal E$.
If $\Cal E$ is ample,then $\det\Cal E$ is ample
and $e=c_1(\Cal E)>0$.
It follows that $d>2$,hence $\Cal E$ is not ample when $d\leq 2$.
On the other hand,$\Cal E$ is ample when $d\geq 5$
by the argument in \cite{F7;(3.13)}.
In general,for any indecomposable vector bundle $\Cal F$
on an elliptic curve,$\Cal F$ is ample if and only if $c_1(\Cal F)>0$
(for a proof,see e.g. \cite{H1}).
Thus when $d=3$ or 4,$\Cal E$ is ample if it is indecomposable.

(3.4)
When $d=3$ or 4,we can find an example of $(M,L)$
by the argument in \cite{F7;(3.12)}.
We can also find an example of $(M,L)$ with $d=6$ as follows.
Let $C$ be a smooth elliptic curve and take a line bundle $\Cal L$
on $C$ with $\deg\Cal L=1$.
We put $\Cal E=\Cal L^{\oplus4}$,then $\Cal E$ is ample,
$c_1(\Cal E)=4,\Bbb P_C(\Cal E)\simeq C\times\Bbb P_\sigma^3$,
and $H(\Cal E)=H_\sigma+\Cal L_{\Bbb P(\Cal E)}$,
where $H_\sigma$ is the pullback of $\Cal O(1)$ on $\Bbb P_\sigma^3$.
Putting $B=-2\Cal L$,we have $\deg B=-2$ and
$2H(\Cal E)+B_{\Bbb P(\Cal E)}=2H_\sigma$.
Then a general member $M$ of $|2H(\Cal E)+B_{\Bbb P(\Cal E)}|$
is smooth and,putting $L=[H(\Cal E)]_M$,
$(M,L)$ becomes an expected example with $d=6$.

(3.5)
{}From now on,we study the case $g(C)=0$.
In this case,$C\simeq\Bbb P_\xi^1$ and we have
$e=d-4$ and $b=8-d$ from the equality in (3.1).
Hence we obtain $d\leq 12$ since $s\geq 0$ and $n\geq 3$.
Furthermore when $d=11$ or 12,we have $n=3$ and when $d=12$,
we have $s=0$ and $\varPhi$ is a $\Bbb P^1\times\Bbb P^1$-bundle
over $\Bbb P_\xi^1$ by \cite{F7;(3.3)}.

(3.6)
We put $P=\Bbb P_C(\Cal E),H=H(\Cal E)$ and $H_\xi=\pi^*\Cal O_C(1)$,
where $\pi$ is the bundle map $\Bbb P_C(\Cal E)\to C$.
Since $\Cal E$ is decomposable,we can describe
$\Cal E\simeq\Cal O(e_0)\oplus\dots\oplus\Cal O(e_n)$,
where $e_0,\dots,e_n\in\Bbb Z,e_0\leq\dots\leq e_n$,and
$\sum_{i=0}^n e_i=e$.
$\Cal O(e_0)\oplus\dots\oplus\Cal O(e_n)$ is denoted by
$\Cal O(e_0,\dots,e_n)$ for simplicity.
We shall classify $\Cal E\simeq\Cal O(e_0,\dots,e_n)$
for each case $d=1,2,\dots,12$.

\proclaim{{\rm(3.7)}Lemma}
$2(e_{n-1}+e_n)<d$ when $e_0\leq 0$.
\endproclaim

\demo{Proof}(cf. \cite{F7;(3.24)}).
A natural surjection $\Cal E\to\Cal O(e_0,\dots,e_{n-1})$
gives a prime divisor $D_1:=\Bbb P(\Cal O(e_0,\dots,e_{n-1}))$ on $P$.
Similarly $\Cal E\to\Cal O(e_0,\dots,e_{n-2},e_n)$
gives a prime divisor $D_2:=\Bbb P(\Cal O(e_0,\dots,e_{n-2},e_n))$ on $P$
and $\Cal E\to\Cal O(e_0,\dots,e_{n-2})$
gives a subvariety $W:=\Bbb P(\Cal O(e_0,\dots,e_{n-2}))$ on $P$.
We have $D_1\in|H-e_nH_\xi|,D_2\in|H-e_{n-1}H_\xi|$,
and $W=D_1\cap D_2$ as schemes.
When $e_0\leq 0$,we have $W\not\subset M$ since $H_W$ is not ample.
Hence $\dim(M\cap W)=n-2$ and
$0<L^{n-2}\{M\cap W\}=H^{n-2}(2H+bH_\xi)(H-e_nH_\xi)(H-e_{n-1}H_\xi)
                   =d-2(e_{n-1}+e_n)$. \qed
\enddemo

(3.8)
Suppose that $d=1$.
We have $e=-3,b=7$,and $M\in|2H+7H_\xi|$.
By (3.7),$\Cal E\simeq\Cal O(-3,0,\dots,0),\Cal O(-2,-1,0,\dots,0)$,
                   or $\Cal O(-1,-1,-1,0,\dots,0)$.

(3.8.1)
When $\Cal E\simeq\Cal O(-1,-1,-1,0,\dots,0)$,we have $n\leq 4$
by the argument in \cite{F7;(3.21)}.
Indeed,we have
  $$ P\simeq \left\{\matrix
    (\xi_0:\xi_1)\times(\sigma_0:\sigma_1:\sigma_2:
                        \sigma_{30}:\sigma_{31}:\dots:
                        \sigma_{n0}:\sigma_{n1})
                 \in\Bbb P_\xi^1\times\Bbb P_\sigma^{2n-2}  \\
    |\xi_0:\xi_1=\sigma_{30}:\sigma_{31}=\dots
                =\sigma_{n0}:\sigma_{n1}
  \endmatrix \right\} $$

\flushpar,$H=H_\sigma-H_\xi$ and $M\in|2H_\sigma+5H_\xi|$.
Thus we can describe

  $$ M=\{q_0(\sigma)\xi_0^5+q_1(\sigma)\xi_0^4\xi_1+\dots
        +q_5(\sigma)\xi_1^5=0~\text{in}~P\} $$

\flushpar,where $q_0,\dots,q_5$ are homogeneous polynomials of degree two in
$\sigma_0,\sigma_1,\dots,\sigma_{n1}$.
In this defining equation of $M$,we put
  $$ \align
    &\sigma_0=a_{00}\xi_0+a_{01}\xi_1,
     \sigma_1=a_{10}\xi_0+a_{11}\xi_1,
     \sigma_2=a_{20}\xi_0+a_{21}\xi_1, \\
    &\sigma_{30}=a_3\xi_0,\sigma_{31}=a_3\xi_1,\dots,
     \sigma_{n0}=a_n\xi_0,\sigma_{n1}=a_n\xi_1
  \endalign $$

\flushpar,where $a_{00},a_{01},\dots,a_n$ are constants.
Then we obtain an equation

$$ Q_0(a)\xi_0^7+Q_1(a)\xi_0^6\xi_1+\dots+Q_7(a)\xi_1^7=0 $$

\flushpar,where $Q_0,\dots,Q_7$ are homogeneous polynomials of degree two in
$(a)=(a_{00},a_{01},\dots,a_n)$.
If $n\geq 5$,then $Q_0(a)=\dots=Q_7(a)=0$ has a non-trivial solution.
We fix such a solution $(a)$ and define a rational map
$\alpha:\Bbb P_\xi^1\to\Bbb P_\sigma^{2n-2}$ by
$$ \align
   \alpha(\xi_0:\xi_1):=(a_{00}\xi_0
   &+a_{01}\xi_1:a_{10}\xi_0+a_{11}\xi_1
                :a_{20}\xi_0+a_{21}\xi_1: \\
   &:a_3\xi_0:a_3\xi_1:\dots:
     a_n\xi_0:a_n\xi_1).
\endalign $$

\flushpar If $\alpha$ is not a morphism,then
$a_{00}:a_{10}:a_{20}=a_{01}:a_{11}:a_{21}$ and $a_3=\dots=a_n=0$.
Since $(a)$ is non-trivial,the equations
$$ \sigma_0:\sigma_1:\sigma_2=a_{00}:a_{10}:a_{20}=a_{01}:a_{11}:a_{21},
   \sigma_{30}=\sigma_{31}=\dots=\sigma_{n0}=\sigma_{n1}=0 $$
determine a point $z$ on $\Bbb P^{2n-2}$.
Let $Z$ be the fiber of a projection
$\Bbb P_\xi^1\times\Bbb P_\sigma^{2n-2}\to
                   \Bbb P_\sigma^{2n-2}$ over $z$.
Then we have $Z\subset M$ by the definition of $Z$,hence
$0<LZ=HZ=(H_\sigma-H_\xi)Z=-1$.
This is a contradiction,thus $\alpha$ is a morphism.
Let $\varGamma$ be the graph of $\alpha$.
Then $\varGamma\subset M$ by the definition of $\alpha$,hence
$0<L\varGamma=H\varGamma=(H_\sigma-H_\xi)\varGamma$.
However,since $H_\sigma\varGamma=H_\xi\varGamma=1$,
this is also a contradiction.
Hence we have proved that $n\leq 4$,thus
$\Cal E\simeq\Cal O(-1,-1,-1,0)$ or $\Cal O(-1,-1,-1,0,0)$.
If $\Cal E\simeq\Cal O(-1,-1,-1,0)$,then
$P\simeq
\{(\xi_0:\xi_1)\times(\sigma_0:\sigma_1:\sigma_2:\sigma_{30}:\sigma_{31})
\in\Bbb P_\xi^1\times\Bbb P_\sigma^4
|\xi_0:\xi_1=\sigma_{30}:\sigma_{31}\}$.
Thus the projection $\mu:P\to\Bbb P_{\sigma}^4$ is the blowing-up of
$\Bbb P_\sigma^4$ with center
$W:=\{\sigma_{30}=\sigma_{31}=0~\text{in}~\Bbb P_\sigma^4\}\simeq\Bbb P^2$.
Since the exceptional divisor $E$ of $\mu$ is a member of $|H_\sigma-H_\xi|$,
we have $M\in|7H_\sigma-5E|$.
Hence $M$ is the strict transform of a hypersurface of degree seven in
$\Bbb P_\sigma^4$,which has singularities with multiplicity five along $W$.

(3.8.2)
When $\Cal E\simeq\Cal O(-2,-1,0,\dots,0)$,we claim that $n\leq 4$.
The following argument is similar to (3.8.1).We have
  $$ P\simeq \left\{\matrix
    (\xi_0:\xi_1)\times(\sigma_0:\sigma_{10}:\sigma_{11}:
                        \sigma_{20}:\sigma_{21}:\sigma_{22}:\dots:
                        \sigma_{n0}:\sigma_{n1}:\sigma_{n2})
                 \in\Bbb P_\xi^1\times\Bbb P_\sigma^{3n-1}  \\
    |\xi_0:\xi_1=\sigma_{10}:\sigma_{11}=\sigma_{20}:\sigma_{21}
                =\sigma_{21}:\sigma_{22}=\dots
                =\sigma_{n0}:\sigma_{n1}=\sigma_{n1}:\sigma_{n2}
  \endmatrix \right\} $$

\flushpar,$H=H_\sigma-2H_\xi$ and $M\in|2H_\sigma+3H_\xi|$.
Thus $M=\{q_0(\sigma)\xi_0^3+q_1(\sigma)\xi_0^2\xi_1+
          q_2(\sigma)\xi_0\xi_1^2+q_3(\sigma)\xi_1^3=0~\text{in}~P\}$
,where $q_0,\dots,q_3$ are quadric polynomials in $(\sigma)$.We put
  $$ \align
    &\sigma_0=a_{00}\xi_0^2+a_{01}\xi_0\xi_1+a_{02}\xi_1^2,
     \sigma_{10}=\xi_0(a_{10}\xi_0+a_{11}\xi_1),
     \sigma_{11}=\xi_1(a_{10}\xi_0+a_{11}\xi_1), \\
    &\sigma_{20}=a_2\xi_0^2,\sigma_{21}=a_2\xi_0\xi_1,\sigma_{22}=a_2\xi_1^2,
     \dots,\sigma_{n0}=a_n\xi_0^2,\sigma_{n1}=a_n\xi_0\xi_1,
     \sigma_{n2}=a_n\xi_1^2
  \endalign $$
Then from the defining equation of $M$ above,we obtain an equation

$$ Q_0(a)\xi_0^7+Q_1(a)\xi_0^6\xi_1+\dots+Q_7(a)\xi_1^7=0 $$

\flushpar,where $Q_0,\dots,Q_7$ are quadric polynomials in
$(a)=(a_{00},a_{01},\dots,a_n)$.
If $n\geq 5$,then $Q_0(a)=\dots=Q_7(a)=0$
has a non-trivial solution $(a)$.
We fix it and define a rational map
$\alpha:\Bbb P_\xi^1\to\Bbb P_\sigma^{3n-1}$ by
$$ \align
   \alpha(\xi_0:\xi_1):=(a_{00}\xi_0^2
   &+a_{01}\xi_0\xi_1+a_{02}\xi_1^2:
   \xi_0(a_{10}\xi_0+a_{11}\xi_1):
   \xi_1(a_{10}\xi_0+a_{11}\xi_1):a_2\xi_0^2  \\
   &:a_2\xi_0\xi_1:a_2\xi_1^2:\dots:
    a_n\xi_0^2:a_n\xi_0\xi_1:a_n\xi_1^2).
\endalign $$

\flushpar If $\alpha$ is not a morphism,then $a_2=\dots=a_n=0$
and for some $(c_0:c_1)\in\Bbb P_\xi^1$,we have
$a_{10}c_0+a_{11}c_1=0$ and $a_{00}c_0^2+a_{01}c_0c_1+a_{02}c_1^2=0$.
In the case $a_{10}=a_{11}=0$,let $Z$ be the fiber of
$\Bbb P_\xi^1\times\Bbb P_\sigma^{3n-1}\to\Bbb P_\sigma^{3n-1}$
over $z:=(1:0:\dots:0)$.
Then we have $Z\subset M$,hence $0<LZ=HZ=(H_\sigma-2H_\xi)Z=-2$.
This is a contradiction,thus $a_{10}\not=0$ or $a_{11}\not=0$.

In this case,$a_{00}\xi_0^2+a_{01}\xi_0\xi_1+a_{02}\xi_1^2$
is devided by $a_{10}\xi_0+a_{11}\xi_1$ in $\Bbb C[\xi_0,\xi_1]$;
we denote by $b_0\xi_0+b_1\xi_1$ its quotient.
We put
$$ Z=\{\sigma_0=b_0\sigma_{10}+b_1\sigma_{11},
       \sigma_{20}=\dots=\sigma_{n2}=0~\text{in}~P\}.
$$
Then $\dim Z=1$ and $Z\subset M$ by the definition of $Z$,hence
$0<LZ=HZ=(H_\sigma-2H_\xi)Z$.
However,since $H_\sigma Z=1$ and $H_\xi Z=1$,this is a contradiction too.
Thus $\alpha$ is a morphism.

Let $\varGamma$ be the graph of $\alpha$.
We have $\varGamma\subset M$ and then
$0<L\varGamma=H\varGamma=(H_\sigma-2H_\xi)\varGamma$.
However,since $H_\sigma\varGamma=2$ and
              $H_\xi\varGamma=1$,this is also a contradiction.
Hence we have proved that $n\leq 4$,thus
$\Cal E\simeq\Cal O(-2,-1,0,0)$ or $\Cal O(-2,-1,0,0,0)$.

(3.8.3)
When $\Cal E\simeq\Cal O(-3,0,\dots,0)$,
we claim that $n\leq 4$ as before.$P$ is isomorphic to
  $$ \left\{\matrix
    (\xi_0:\xi_1)\times(\sigma_0:\sigma_{10}:\sigma_{11}:
                        \sigma_{12}:\sigma_{13}:\dots:
                        \sigma_{n0}:\sigma_{n1}:\sigma_{n2}:\sigma_{n3})
                 \in\Bbb P_\xi^1\times\Bbb P_\sigma^{4n}  \\
    |\xi_0:\xi_1=\sigma_{10}:\sigma_{11}=\sigma_{11}:\sigma_{12}
                =\sigma_{12}:\sigma_{13}=\dots
                =\sigma_{n0}:\sigma_{n1}=\sigma_{n1}:\sigma_{n2}
                =\sigma_{n2}:\sigma_{n3}
  \endmatrix \right\} $$
,$H=H_\sigma-3H_\xi$ and $M\in|2H_\sigma+H_\xi|$.
Thus $M=\{q_0(\sigma)\xi_0+q_1(\sigma)\xi_1=0~\text{in}~P\}$,
where $q_0$ and $q_1$ are quadric polynomials in $(\sigma)$.We put
  $$ \align
    &\sigma_0=a_{00}\xi_0^3+a_{01}\xi_0^2\xi_1
                           +a_{02}\xi_0\xi_1^2
                           +a_{03}\xi_1^3, \\
    &\sigma_{10}=a_1\xi_0^3,\sigma_{11}=a_1\xi_0^2\xi_1,
     \sigma_{12}=a_1\xi_0\xi_1^2,\sigma_{13}=a_1\xi_1^3,
     \dots, \\
    &\sigma_{n0}=a_n\xi_0^3,\sigma_{n1}=a_n\xi_0^2\xi_1,
     \sigma_{n2}=a_n\xi_0\xi_1^2,\sigma_{n3}=a_n\xi_1^3.
  \endalign $$
Then from the defining equation of $M$ above,we obtain an equation
  $$ Q_0(a)\xi_0^7+Q_1(a)\xi_0^6\xi_1+\dots+Q_7(a)\xi_1^7=0 $$
,where $Q_0,\dots,Q_7$ are quadric polynomials in
$(a)=(a_{00},a_{01},\dots,a_n)$.
If $n\geq 5$,then $Q_0(a)=\dots=Q_7(a)=0$ has a non-trivial solution $(a)$.
We fix it and define a rational map
$\alpha:\Bbb P_\xi^1\to\Bbb P_\sigma^{4n}$ by
  $$ \align
     \alpha(\xi_0:\xi_1):=(a_{00}\xi_0^3
     &+a_{01}\xi_0^2\xi_1+a_{02}\xi_0\xi_1^2+a_{03}\xi_1^3:
      a_1\xi_0^3:a_1\xi_0^2\xi_1:a_1\xi_0\xi_1^2:a_1\xi_1^3 \\
     &:\dots:a_n\xi_0^3:a_n\xi_0^2\xi_1:a_n\xi_0\xi_1^2:a_n\xi_1^3).
     \endalign $$
If $\alpha$ is not a morphism,then $a_1=\dots=a_n=0$.
Let $Z$ be the fiber of
$\Bbb P_\xi^1\times\Bbb P_\sigma^{4n}\to\Bbb P_\sigma^{4n}$
over $z:=(1:0:\dots:0)$.
We have $Z\subset M$ and then $0<LZ=HZ=(H_\sigma-3H_\xi)Z=-3$.
This is a contradiction,hence $\alpha$ is a morphism.
Let $\varGamma$ be the graph of $\alpha$.
We have $\varGamma\subset M$ and then
$0<L\varGamma=H\varGamma=(H_\sigma-3H_\xi)\varGamma$.
However,since $H_\sigma\varGamma=3$ and
              $H_\xi\varGamma=1$,this is also a contradiction.
Hence we have proved that $n\leq 4$,thus
$\Cal E\simeq\Cal O(-3,0,0,0)$ or $\Cal O(-3,0,0,0,0)$.

(3.9)
Now we study the case $d=2$.
We have $e=-2,b=6$,and $M\in|2H+6H_\xi|$.
By (3.7),$\Cal E\simeq\Cal O(-2,0,\dots,0)$ or $\Cal O(-1,-1,0,\dots,0)$.

(3.9.1)
When $\Cal E\simeq\Cal O(-1,-1,0,\dots,0)$,
we have $n\leq 4$ as in (3.8.1).
Hence $\Cal E\simeq\Cal O(-1,-1,0,0)$ or $\Cal O(-1,-1,0,0,0)$.

(3.9.2)
When $\Cal E\simeq\Cal O(-2,0,\dots,0)$,we have $n\leq 4$ as in (3.8.2).
Hence $\Cal E\simeq\Cal O(-2,0,0,0)$ or $\Cal O(-2,0,0,0,0)$.

(3.10)
Suppose that $d=3$.
Then $e=-1,b=5$,and $M\in|2H+5H_\xi|$.
{}From (3.7),we have $\Cal E\simeq\Cal O(-2,0,\dots,0,1)$,
$\Cal E\simeq\Cal O(-1,-1,0,\dots,0,1)$,
or $\Cal E\simeq\Cal O(-1,0,\dots,0)$.

(3.10.1)
When $\Cal E\simeq\Cal O(-1,0,\dots,0)$,
we have $n\leq 4$ as in (3.8.1).
Hence $\Cal E\simeq\Cal O(-1,0,0,0)$ or $\Cal O(-1,0,0,0,0)$.

(3.10.2)
When $\Cal E\simeq\Cal O(-1,-1,0,\dots,0,1)$,
we have $n\leq 4$ by the argument in \cite{F7;(3.23.2)}
which is similar to (3.8.1).
Hence $\Cal E\simeq\Cal O(-1,-1,0,1)$ or $\Cal O(-1,-1,0,0,1)$.

(3.10.3)
When $\Cal E\simeq\Cal O(-2,0,\dots,0,1)$,
we have $n\leq 4$ as in (3.8.2) and (3.10.2).
Hence $\Cal E\simeq\Cal O(-2,0,0,1)$ or $\Cal O(-2,0,0,0,1)$.

The next lemma is useful for $d\geq 4$.

\proclaim{{\rm(3.11)}Lemma}
When $d\geq 4$,$-1$ does not appear twice in $\{e_0,\dots,e_n\}$.
\endproclaim

We can prove this lemma by the argument in \cite{F7;(3.18)}.

(3.12)
Now we study the case $d=4$.
We have $e=0,b=4$,and $M\in|2H+4H_\xi|$.
By (3.7) and (3.11),
$\Cal E\simeq\Cal O(-1,0,\dots,0,1)$ or $\Cal O(0,\dots,0)$.

(3.12.1)
When $\Cal E\simeq\Cal O(-1,0,\dots,0,1)$,
we have $n\leq 4$ as in (3.10.2).
Hence $\Cal E\simeq\Cal O(-1,0,0,1)$ or $\Cal O(-1,0,0,0,1)$.

(3.12.2)
When $\Cal E\simeq\Cal O(0,\dots,0)$,
by the argument in \cite{F7;(3.23.1)},we have $n\leq 4,
P\simeq\Bbb P_\xi^1\times\Bbb P_\sigma^n$,Bs$|L|=\phi$,
and the morphism $\varphi:M\to\Bbb P_\sigma^n$ defined by $|L|$
is a finite morphism of degree four.
Conversely,any general member $M$ of $|2H_\sigma+4H_\xi|$ on $P$
does not contain any fiber of the projection $P\to\Bbb P_\sigma^n$,
thus $L:=H_M$ is ample and $(M,L)$ is a polarized manifold
of the above type.

The next lemma is useful for $d\geq 5$.

\proclaim{{\rm(3.13)}Lemma}
$e_0\geq -1$ when $d\geq 5$.
\endproclaim

We can prove this lemma by the argument in \cite{F7;(3.19)}.

Similarly we obtain the following two lemmas.

\proclaim{{\rm(3.14)}Lemma}
$e_0\geq 0$ when $d\geq 7$.
\endproclaim

\proclaim{{\rm(3.15)}Lemma}
$e_0\geq 1$ when $d\geq 9$.
\endproclaim

(3.16)
Now we study the case $d=5$.
We have $e=1,b=3$,and $M\in|2H+3H_\xi|$.
By (3.11) and (3.13),$\Cal E\simeq\Cal O(-1,0,\dots,0,2)$,
$\Cal O(-1,0,\dots,0,1,1)$,or $\Cal O(0,\dots,0,1)$.

(3.16.1)
When $\Cal E\simeq\Cal O(-1,0,\dots,0,2)$,
we have $n\leq 3$ as in (3.10.2),
hence $\Cal E\simeq\mathbreak\Cal O(-1,0,0,2)$.
Furthermore Bs$|L|$ is one point as in \cite{F7;(3.23.2)}.

(3.16.2)
When $\Cal E\simeq\Cal O(-1,0,\dots,0,1,1)$,
we have $n\leq 4$ and Bs$|L|$ is one point as in (3.16.1).
Thus $\Cal E\simeq\Cal O(-1,0,1,1)$ or $\Cal O(-1,0,0,1,1)$.

(3.16.3)
When $\Cal E\simeq\Cal O(0,\dots,0,1)$,
by the argument in \cite{F7;(3.24)},
we have $n\leq 4$ and $|L|$ makes $M$ the normalization of
a hypersurface of degree five in $\Bbb P^{n+1}$,
which has triple points along a $\Bbb P^2$ in $\Bbb P^{n+1}$.

(3.17)
Suppose that $d=6$.
We have $e=2,b=2$,and $M\in|2H+2H_\xi|$.
By (3.7),(3.11),and (3.13),$\Cal E\simeq\Cal O(-1,0,\dots,0,1,1,1)$,
$\Cal O(0,\dots,0,1,1)$,or $\Cal O(0,\dots,0,2)$.

(3.17.1)
When $\Cal E\simeq\Cal O(-1,0,\dots,0,1,1,1)$,
we show that $n=3$ similarly as in (3.7).
Natural surjections $\Cal E\to\Cal O(e_0,\dots,e_{n-1})$,
                    $\Cal E\to\Cal O(e_0,\dots,e_{n-2},e_n)$,and
                    $\Cal E\to\Cal O(e_0,\dots,\mathbreak
                                     e_{n-3},e_{n-1},e_n)$
give prime divisors $D_1:=\Bbb P(\Cal O(e_0,\dots,e_{n-1}))$,
                    $D_2:=\Bbb P(\Cal O(e_0,\dots,e_{n-2},e_n))$
                        and
                    $D_3:=\Bbb P(\Cal O(e_0,\dots,e_{n-3},e_{n-1},e_n))$
respectively.
A natural surjection $\Cal E\to\Cal O(e_0,\dots,\mathbreak
                                      e_{n-3})$
gives a subvariety $W:=\Bbb(\Cal O(e_0,\dots,e_{n-3}))$
of $P=\Bbb P(\Cal E)$.
We have $D_1\in|H-e_nH_\xi|$,
        $D_2\in|H-e_{n-1}H_\xi|$,
        $D_3\in|H-e_{n-2}H_\xi|$,
and $W=D_1\cap D_2\cap D_3$ as schemes.
Since $H_W$ is not ample,we have $W\not\subset M$,hence
$\dim(M\cap W)=n-3$ and
$0<L^{n-3}\{M\cap W\}=H^{n-3}(2H+2H_\xi)(H-H_\xi)^3=2e-4=0$
if $n\leq 4$.
This is a contradiction,thus we have $n=3$ and
$\Cal E\simeq\Cal O(-1,1,1,1)$.
By the argument in \cite{F7;(3.26)},$M$ is a double covering of
$\Bbb P_\xi^1\times\Bbb P_\sigma^2$
and its branch locus is a smooth member of $|4H_\xi+2H_\sigma|$.
We have also $L=[H_\xi+H_\sigma]_M$.

(3.17.2)
When $\Cal E\simeq\Cal O(0,\dots,0,1,1)$,
we have $n\leq 4$ as in (3.16.3),hence
$\Cal E\simeq\Cal O(0,0,1,1)$ or $\Cal O(0,0,0,1,1)$.
We show the existence of $(M,L)$.
When $\Cal E\simeq\Cal O(0,0,1,1)$,we have
$P\simeq\{(\xi_0:\xi_1)\times(\sigma_0:\sigma_1:
                              \sigma_{20}:\sigma_{21}:
                              \sigma_{30}:\sigma_{31})
\in\Bbb P_\xi^1\times\Bbb P_\sigma^5
|\xi_0:\xi_1=\sigma_{20}:\sigma_{21}
            =\sigma_{30}:\sigma_{31}\}$
and $H=H_\sigma$.
Let $M$ be a general member of $|2H_\sigma+2H_\xi|$
and put $L=[H_\sigma]_M$.
Then Bs$|L|=\phi$ and the restriction of the projection
$P\to\Bbb P_\sigma^5$ to $M$ is the morphism
$\varphi$ defined by $|L|$.
If $\varphi:M\to\varphi(M)$ is not finite,
$M$ contains a fiber $Z$ of the projection
$P\to\Bbb P_\sigma^5$ over one point $z$ on the line
$l:=\{\sigma_{20}=\sigma_{21}=\sigma_{30}=\sigma_{31}
                 =0~\text{in}~\Bbb P_\sigma^5\}$
Using homogeneous polynomials $q_0,q_1$,and $q_2$
of degree two in $(\sigma)$,we can describe that
$M=\{q_0(\sigma)\xi_0^2+q_1(\sigma)\xi_0\xi_1+
     q_2(\sigma)\xi_1^2=0~\text{in}~P\}$.
Then $Z\subset M$ if and only if $q_0(z)=q_1(z)=q_2(z)=0$.
Thus if we choose $q_0,q_1$,and $q_2$ generally to satisfy that
$l\cap\{q_0(\sigma)=q_1(\sigma)=q_2(\sigma)
                   =0~\text{in}~\Bbb P_\sigma^5\}=\phi$,
then $\varphi$ becomes finite and $L$ is ample.
Similarly we can find examples of $(M,L)$
such that $\Cal E\simeq\Cal O(0,0,0,1,1)$.

(3.17.3)
When $\Cal E\simeq\Cal O(0,\dots,0,2)$,
we have $n\leq 3$ as in (3.16.3),hence
$\Cal E\simeq\Cal O(0,0,0,2)$.
We can show the existence of $(M,L)$ similarly as above.

When $d\geq 7$,the situation is much simpler.

\proclaim{{\rm(3.18)}Lemma}
 Bs$|L|=\phi$ and $L$ is very ample when $d\geq 7$.
\endproclaim

We can prove this lemma similarly as in \cite{F7;(3.31)}.
This lemma tells us that our results overlap \cite{I;Theorem 4.3},
but our method is different from his.

(3,19)
Now we study the case $d=7$.
We have $e=3,b=1$,and $M\in|2H+H_\xi|$.
Furthermore $e_0\geq 0$ by (3.14) and
            $e_2\geq 1$ by the argument in \cite{F7;(3.25)}.
Hence $\Cal E\simeq\Cal O(0,0,1,2),\Cal O(0,1,1,1)$,or $\Cal O(0,0,1,1,1)$.
In each case,$(M,L)$ exists similarly as in (3.17.2).
By the morphism defined by $|L|$,$M$ is isomorphic to a manifold
of degree seven in $\Bbb P^{n+3}$.

(3.20)
Suppose that $d=8$.
We have $e=4,b=0$,and $M\in|2H|$.
Furthermore $e_0\geq 0$ by (3.14) and
            $e_1\geq 1$ by the argument in \cite{F7;(3.26)}.
Hence $\Cal E\simeq\Cal O(0,1,1,2),\Cal O(0,1,1,1,1)$,or $\Cal O(1,1,1,1)$.

(3.20.1)
When $\Cal E\simeq\Cal O(1,1,1,1)$,we have
$P\simeq\Bbb P_\xi^1\times\Bbb P_\sigma^3,H=H_\xi+H_\sigma$,
and $M\in|2H_\sigma+2H_\xi|$.
Hence $M$ is a smooth divisor of bidegree (2,2) on $P$.
Conversely,let $M$ be a general member of $|2H_\xi+2H_\sigma|$
and put $L=[H_\xi+H_\sigma]_M$.
Since $\Cal E$ is ample,$L$ is ample and
$(M,L)$ is a polarized manifold of the above type.

(3.20.2)
When $\Cal E\simeq\Cal O(0,1,1,1,1)$,
by the argument in \cite{F7;(3.26)},
$M$ is a double covering of $\Bbb P_\xi^1\times\Bbb P_\sigma^3$
and its branch locus is a smooth member of $|2H_\xi+2H_\sigma|$.
We have also $L=[H_\xi+H_\sigma]_M$.

(3.20.3)
Even when $\Cal E\simeq\Cal O(0,1,1,2)$,
by the argument in \cite{F7;(3.26)},
we have a morphism $h:M\to\Bbb P_\xi^1\times\Bbb P_\sigma^3$
and $L=h^*(H_\xi+H_\sigma)$.
Since $L$ is ample,$h:M\to h(M)$ is finite and
$h(M)\in|a_1H_\xi+a_2H_\sigma|$ for some non-negative integers
$a_1$ and $a_2$.
Then $8=L^3=(\deg h)\cdot[H_\xi+H_\sigma]_{h(M)}^3
           =(\deg h)(a_1+3a_2)$.
{}From the construction of $h$,
we get $\deg h=2$ and $a_1=a_2=1$.
Hence $h(M)\in|H_\xi+H_\sigma|$ and $M\to h(M)$ is a double covering.

(3.21)
Suppose that $d=9$.
We have $e=5,b=-1$,and $M\in|2H-H_\xi|$.
Since $e_0\geq 1$ by (3.15),
$\Cal E\simeq\Cal O(1,1,1,2)$ or $\Cal O(1,1,1,1,1)$.

(3.21.1)
When $\Cal E\simeq\Cal O(1,1,1,1,1)$,similarly as in \cite{F7;(3.27)},
the restriction of the projection
$P\simeq\Bbb P_\xi^1\times\Bbb P_\sigma^4\to\Bbb P_\sigma^4$
to $M$ is a blowing-up of $\Bbb P_\sigma^4$
and its center is a complete intersection of two hyperquadrics
in $\Bbb P_\sigma^4$.

(3.21.2)
When $\Cal E\simeq\Cal O(1,1,1,2)$,we have
$P\simeq\{(\xi_0:\xi_1)\times(\sigma_0:\sigma_1:\sigma_2:
                              \sigma_{30}:\sigma_{31})
  \in\Bbb P_\xi^1\times\Bbb P_\sigma^4
  |\xi_0:\xi_1=\sigma_{30}:\sigma_{31}\}$,
hence $P$ is the blowing-up of $\Bbb P_\sigma^4$ with center
$\{\sigma_{30}=\sigma_{31}=0~\text{in}~\Bbb P_\sigma^4\}$.
The exceptional divisor $E$ is
$\{\sigma_{30}=\sigma_{31}=0~\text{in}~P\}\in|H_\sigma-H_\xi|$,
thus $M\in|3H_\sigma-E|$ and $M$ is the strict transform of
a smooth hypercubic in $\Bbb P_\sigma^4$.

(3.22)
Suppose that $d=10$.
We have $e=6,b=-2$,and $M\in|2H-2H_\xi|$.
Since $e_0\geq 1$ by (3.15),
$\Cal E\simeq\Cal O(1,1,1,3),\Cal O(1,1,2,2),\Cal O(1,1,1,1,2)$,
or $\Cal O(1,1,1,1,1,1)$.

(3.22.1)
When $\Cal E\simeq\Cal O(1,1,1,1,1,1)$,we have
$P\simeq\Bbb P_\xi^1\times\Bbb P_\sigma^5,H=H_\xi+H_\sigma,
 M\in|2H_\sigma|$,and $L=[H_\xi+H_\sigma]_M$.
Hence $M\simeq\Bbb P_\xi^1\times Q$,where $Q$ is a smooth
hyperquadric in $\Bbb P_\sigma^5$.

(3.22.2)
When $\Cal E\simeq\Cal O(1,1,1,1,2)$,by the argument in \cite{F7;(3.28)},
we have $M$ is the blowing-up of a hyperquadric in $\Bbb P_\sigma^5$
and its center is a smooth quadric surface.

(3.22.3)
When $\Cal E\simeq\Cal O(1,1,2,2)$,we have
$P\simeq\{(\xi_0:\xi_1)\times(\sigma_0:\sigma_1:
                              \sigma_{20}:\sigma_{21}:
                              \sigma_{30}:\sigma_{31})
  \in\Bbb P_\xi^1\times\Bbb P_\sigma^5
  |\xi_0:\xi_1=\sigma_{20}:\sigma_{21}=\sigma_{30}:\sigma_{31}\}$,
$H=H_\xi+H_\sigma,M\in|2H_\sigma|$,and $L=[H_\xi+H_\sigma]_M$.
Since $\Cal E$ is ample,$H$ is ample and then
$L$ is ample for any general member $M$ of $|2H_\sigma|$.
Because of (3.18),$M$ is embedded in $\Bbb P^9$ as a manifold
of degree nine by the morphism defined by $|L|$.
On the other hand,the restriction of the projection
$\mu:P\to\Bbb P_\sigma^5$ to $M$ is the morphism defined by $|L-H_\xi|$,
and $M$ is birationally mapped onto $\mu(M)$.
We have
$10=L^3=3[H_\xi]_M[H_\sigma]_M^2+[H_\sigma]_M^3$
and $[H_\xi]_M[H_\sigma]_M^2=2$
since $M\to\Bbb P_\xi^1$ is a hyperquadric fibration.
Thus the degree of $\mu(M)$ is four.
Furthermore,since $\mu(P)=\{\sigma_{20}\sigma_{31}-\sigma_{30}\sigma_{21}=0
                            ~\text{in}~\Bbb P_\sigma^5\}$
and $M\in|2H_\sigma|$,
$\mu(M)$ is a complete intersection of two hyperquadrics
in $\Bbb P_\sigma^5$.
Even when $\Cal E\simeq\Cal O(1,1,1,3)$,
we have the same result as above.

(3.23)
Suppose that $d=11$.
We have $e=7,b=-3$,and $M\in|2H-3H_\xi|$.
Since $e_0\geq 1$ by (3.15),and since $n=3$ by (3.5),
$\Cal E\simeq\Cal O(1,1,1,4),\Cal O(1,1,2,3)$ or $\Cal O(1,2,2,2)$.

(3.23.1)
When $\Cal E\simeq\Cal O(1,1,1,4)$,we claim that $(M,L)$ does not exist.
Assume that $(M,L)$ exists.
A natural surjection $\Cal E\to\Cal O(1,1,1)$ gives a prime divisor
$W:=\Bbb P(\Cal O(1,1,1))$ on $P$.
Since $W\simeq\Bbb P_\xi^1\times\Bbb P_\sigma^2$,
$H_W=H_\xi+H_\sigma$,and $W\not\subset M$,
we have $[M]_W=M\cap W\in|2H_W-3H_\xi|=|2H_\sigma-H_\xi|$.
This is a contradiction,thus we have proved the claim.

(3.23.2)
Even when $\Cal E\simeq\Cal O(1,1,2,3)$,
we can show that $(M,L)$ does not exist. We have
$P\simeq\{(\xi_0:\xi_1)\times(\sigma_0:\sigma_1:
                              \sigma_{20}:\sigma_{21}:
                              \sigma_{30}:\sigma_{31}:\sigma_{32})
  \in\Bbb P_\xi^1\times\Bbb P_\sigma^6
 |\xi_0:\xi_1=\sigma_{20}:\sigma_{21}
             =\sigma_{30}:\sigma_{31}
             =\sigma_{31}:\sigma_{32}\}$
and $H=H_\sigma+H_\xi$.
Assume that there exists a smooth member $M$ of $|2H_\sigma-H_\xi|$.
Then there is an exact sequence of normal bundles
  $$ 0\to\Cal N_{B/M}\to\Cal N_{B/P}\to [\Cal N_{M/P}]_B\to 0 $$
,where $B:=\text{Bs}|2H_\sigma-H_\xi|
          =\{\sigma_{20}=\sigma_{21}=\sigma_{30}
          =\sigma_{31}=\sigma_{32}=0~\text{in}~P\}
          \simeq\Bbb P(\Cal O(1,1))$.
Since $B$ is the complete intersection of
$D_1:=\{\sigma_{20}=\sigma_{21}=0~\text{in}~P\}\simeq\Bbb P(\Cal O(1,1,3))$
and
$D_2:=\{\sigma_{30}=\sigma_{31}=\sigma_{32}=0~\text{in}~P\}
                                               \simeq\Bbb P(\Cal O(1,1,2))$
,we have $\Cal N_{B/P}\simeq [\Cal N_{D_1/P}]_B\oplus [\Cal N_{D_2/P}]_B
                      \simeq [H_\sigma-H_\xi]_B\oplus [H_\sigma-2H_\xi]_B$.
Also we have $\Cal N_{M/P}\simeq [2H_\sigma-H_\xi]_B$.
Then the morphism
$\varphi:[H_\sigma-H_\xi]_B\oplus [H_\sigma-2H_\xi]_B
               \to [2H_\sigma-H_\xi]_B$
corresponding to
$\Cal N_{B/P}\to [\Cal N_{M/P}]_B$
is given by some $\varphi_1\in H^0(B,[H_\sigma]_B)$ and
                 $\varphi_2\in H^0(B,[H_\sigma+H_\xi]_B)$.
Since $[H_\sigma]_B[H_\sigma+H_\xi]_B=1$,
$\varphi_1$ and $\varphi_2$ have a common zero point,at which
$\varphi$ is not surjective.
This yields a contradiction and $(M,L)$ does not exist.

(3.23.3)
When $\Cal E\simeq\Cal O(1,2,2,2)$,we can show the existence of $(M,L)$.
We have
$P\simeq\{(\xi_0:\xi_1)\times(\sigma_0:
                              \sigma_{10}:\sigma_{11}:
                              \sigma_{20}:\sigma_{21}:
                              \sigma_{30}:\sigma_{31})
  \in\Bbb P_\xi^1\times\Bbb P_\sigma^6
 |\xi_0:\xi_1=\sigma_{10}:\sigma_{11}
             =\sigma_{20}:\sigma_{21}
             =\sigma_{30}:\sigma_{31}\}$
and $H=H_\sigma+H_\xi$.
Putting $U_i=\{\xi_i\not=0~\text{in}~P\}$ and
        $V_j=\{\sigma_j\not=0~\text{in}~P\}$,
we take a rational section
$s_1:=\{(U_i\cap V_j,\frac{\sigma_0^2}{\xi_0}\cdot
                   \frac{\xi_i}{\sigma_j^2})\}_{i,j}$
of $2H_\sigma-H_\xi$.
Note that
$h^0(P,2H-3H_\xi)=h^0(\Bbb P_\xi^1,S^2(\Cal E)\otimes[-3H_\xi])=15$.
Let $f_1,\dots,f_{15}$ be rational functions on $P$ such that
  $$ \alignat 2
     f_1&=\frac{\xi_0}{\sigma_0^2}\cdot
          \frac{\sigma_0\sigma_{10}}{\xi_0}
         =\frac{\xi_0}{\sigma_0^2}\cdot
          \frac{\sigma_0\sigma_{11}}{\xi_1},&~
     f_2&=\frac{\xi_0}{\sigma_0^2}\cdot
          \frac{\sigma_0\sigma_{20}}{\xi_0}
         =\frac{\xi_0}{\sigma_0^2}\cdot
          \frac{\sigma_0\sigma_{21}}{\xi_1},\\
     f_3&=\frac{\xi_0}{\sigma_0^2}\cdot
          \frac{\sigma_0\sigma_{30}}{\xi_0}
         =\frac{\xi_0}{\sigma_0^2}\cdot
          \frac{\sigma_0\sigma_{31}}{\xi_1},&~
     f_4&=\frac{\xi_0}{\sigma_0^2}\cdot
          \frac{\sigma_{10}^2}{\xi_0}
         =\frac{\xi_0}{\sigma_0^2}\cdot
          \frac{\sigma_{10}\sigma_{11}}{\xi_1},\\
     f_5&=\frac{\xi_0}{\sigma_0^2}\cdot
          \frac{\sigma_{10}\sigma_{11}}{\xi_0}
         =\frac{\xi_0}{\sigma_0^2}\cdot
          \frac{\sigma_{11}^2}{\xi_1},&~
     f_6&=\frac{\xi_0}{\sigma_0^2}\cdot
          \frac{\sigma_{10}\sigma_{20}}{\xi_0}
         =\frac{\xi_0}{\sigma_0^2}\cdot
          \frac{\sigma_{10}\sigma_{21}}{\xi_1},\\
     f_7&=\frac{\xi_0}{\sigma_0^2}\cdot
          \frac{\sigma_{10}\sigma_{21}}{\xi_0}
         =\frac{\xi_0}{\sigma_0^2}\cdot
          \frac{\sigma_{11}\sigma_{21}}{\xi_1},&~
     f_8&=\frac{\xi_0}{\sigma_0^2}\cdot
          \frac{\sigma_{10}\sigma_{30}}{\xi_0}
         =\frac{\xi_0}{\sigma_0^2}\cdot
          \frac{\sigma_{10}\sigma_{31}}{\xi_1},\\
     f_9&=\frac{\xi_0}{\sigma_0^2}\cdot
          \frac{\sigma_{10}\sigma_{31}}{\xi_0}
         =\frac{\xi_0}{\sigma_0^2}\cdot
          \frac{\sigma_{11}\sigma_{31}}{\xi_1},&~
  f_{10}&=\frac{\xi_0}{\sigma_0^2}\cdot
          \frac{\sigma_{20}^2}{\xi_0}
         =\frac{\xi_0}{\sigma_0^2}\cdot
          \frac{\sigma_{20}\sigma_{21}}{\xi_1},\\
  f_{11}&=\frac{\xi_0}{\sigma_0^2}\cdot
          \frac{\sigma_{20}\sigma_{21}}{\xi_0}
         =\frac{\xi_0}{\sigma_0^2}\cdot
          \frac{\sigma_{21}^2}{\xi_1},&~
  f_{12}&=\frac{\xi_0}{\sigma_0^2}\cdot
          \frac{\sigma_{20}\sigma_{30}}{\xi_0}
         =\frac{\xi_0}{\sigma_0^2}\cdot
          \frac{\sigma_{20}\sigma_{31}}{\xi_1},\\
  f_{13}&=\frac{\xi_0}{\sigma_0^2}\cdot
          \frac{\sigma_{20}\sigma_{31}}{\xi_0}
         =\frac{\xi_0}{\sigma_0^2}\cdot
          \frac{\sigma_{21}\sigma_{31}}{\xi_1},&~
  f_{14}&=\frac{\xi_0}{\sigma_0^2}\cdot
          \frac{\sigma_{30}^2}{\xi_0}
         =\frac{\xi_0}{\sigma_0^2}\cdot
          \frac{\sigma_{30}\sigma_{31}}{\xi_1},\\
  f_{15}&=\frac{\xi_0}{\sigma_0^2}\cdot
          \frac{\sigma_{30}\sigma_{31}}{\xi_0}
         =\frac{\xi_0}{\sigma_0^2}\cdot
          \frac{\sigma_{31}^2}{\xi_1}.
\endalignat $$
Then $\Bbb C\langle f_1,\dots,f_{15}\rangle$,
the vector space spanned by $f_1,\dots,f_{15}$ over $\Bbb C$,
is isomorphic to $H^0(P,2H_\sigma-H_\xi)$
by mapping each $f_i$ to $f_i\cdot s_1$.
Thus we can describe
  $$ |2H_\sigma-H_\xi|=\{~\text{div}(f\cdot s_1)|
     f\in\Bbb C\langle f_1,\dots,f_{15}\rangle-0\} $$
,where div$(f\cdot s_1)$ is an effective divisor defined by
a regular section $f\cdot s_1$ of $2H_\sigma-H_\xi$.
Since
     $\text{Bs}|2H_\sigma-H_\xi|=\{\sigma_{10}=\sigma_{11}=\dots
                                              =\sigma_{31}=0~\text{in}~P\}
      \simeq\Bbb P_\xi^1\times\{(1:0:\dots:0)\}$,
if we take $f=\sum_{i=1}^{15} c_if_i\in\Bbb C
                                       \langle f_1,\dots,f_{15}\rangle$
with $(c_1,c_2,c_3)\not=(0,0,0)$,
$\text{div}(f\cdot s_1)$ is nonsingular along
$\text{Bs}|2H_\sigma-H_\xi|$.
Thus any general member $M$ of $|2H_\sigma-H_\xi|$ is smooth
by Bertini's theorem.
For such $M$,$L:=H_M$ is ample since $\Cal E$ is ample,
hence $(M,L)$ is a polarized manifold as desired.
Furthermore,similarly as in (3.16.3),
$|L-H_\xi|$ makes $M$ a desingularization of a variety
of degree five in $\Bbb P_\sigma^6$.

(3.24)
Suppose that $d=12$.
We have $e=8,b=-4$,and $M\in|2H-4H_\xi|$.
Since $e_0\geq 1$ by (3.15),and since $n=3$ by (3.5),
$\Cal E\simeq\Cal O(1,1,1,5),\Cal O(1,1,2,4),\Cal O(1,1,3,3),
             \Cal O(1,2,2,3)$,or $\Cal O(2,2,2,2)$.

(3.24.1)
When $\Cal E\simeq\Cal O(2,2,2,2)$,we have
$P\simeq\Bbb P_\xi^1\times\Bbb P_\sigma^3,
 H=H_\sigma+2H_\xi,M\in|2H_\sigma|$,
and $L=[H_\sigma+H_\xi]_M$.
Hence $M\simeq\Bbb P_\xi^1\times Q$,where $Q$ is
a smooth quadric surface in $\Bbb P_\sigma^3$.
Since $Q\simeq\Bbb P_\mu^1\times\Bbb P_\lambda^1$,
we have $M\simeq\Bbb P_\xi^1\times
                \Bbb P_\mu^1\times
                \Bbb P_\lambda^1$
and $L=2H_\xi+H_\mu+H_\lambda$.

(3.24.2)
When $\Cal E\simeq\Cal O(1,1,1,5)$,$(M,L)$ does not exist
by the argument in (3.23.1).

(3.24.3)
Even when $\Cal E\simeq\Cal O(1,1,2,4)$,
we can show that $(M,L)$ does not exist similarly as in (3.23.2).

(3.24.4)
When $\Cal E\simeq\Cal O(1,2,2,3)$,we can show the existence of $(M,L)$
similarly as in (3.23.3).
In fact,we have
$P\simeq\{(\xi_0:\xi_1)\times(\sigma_0:
                              \sigma_{10}:\sigma_{11}:
                              \sigma_{20}:\sigma_{21}:
                              \sigma_{30}:\sigma_{31}:\sigma_{32})
  \in\Bbb P_\xi^1\times\Bbb P_\sigma^7
 |\xi_0:\xi_1=\sigma_{10}:\sigma_{11}
             =\sigma_{20}:\sigma_{21}
             =\sigma_{30}:\sigma_{31}
             =\sigma_{31}:\sigma_{32}\}$,
$H=H_\sigma+H_\xi$,and
$h^0(P,2H-4H_\xi)=h^0(\Bbb P_\xi^1,S^2(\Cal E)\otimes[-4H_\xi])=11$.
We take a rational section
$s_2:=\{(U_i\cap V_j,\frac{\sigma_0^2}{\xi_0^2}\cdot
                  \frac{\xi_i^2}{\sigma_j^2})\}_{i,j}$
of $2H_\sigma-2H_\xi$,
where $U_i$ and $V_j$ are the same as in (3.23.3).
Let $f_1,\dots,f_{11}$ be rational functions on $P$ such that
  $$ \alignat 4
     f_1&=\frac{\xi_0^2}{\sigma_0^2}\cdot
          \frac{\sigma_0\sigma_{30}}{\xi_0^2},&~
     f_2&=\frac{\xi_0^2}{\sigma_0^2}\cdot
          \frac{\sigma_{10}^2}{\xi_0^2},&~
     f_3&=\frac{\xi_0^2}{\sigma_0^2}\cdot
          \frac{\sigma_{10}\sigma_{20}}{\xi_0^2},&~
     f_4&=\frac{\xi_0^2}{\sigma_0^2}\cdot
          \frac{\sigma_{10}\sigma_{30}}{\xi_0^2},\\
     f_5&=\frac{\xi_0^2}{\sigma_0^2}\cdot
          \frac{\sigma_{10}\sigma_{31}}{\xi_0^2},&~
     f_6&=\frac{\xi_0^2}{\sigma_0^2}\cdot
          \frac{\sigma_{20}^2}{\xi_0^2},&~
     f_7&=\frac{\xi_0^2}{\sigma_0^2}\cdot
          \frac{\sigma_{20}\sigma_{30}}{\xi_0^2},&~
     f_8&=\frac{\xi_0^2}{\sigma_0^2}\cdot
          \frac{\sigma_{20}\sigma_{31}}{\xi_0^2},\\
     f_9&=\frac{\xi_0^2}{\sigma_0^2}\cdot
          \frac{\sigma_{30}^2}{\xi_0^2},&~
  f_{10}&=\frac{\xi_0^2}{\sigma_0^2}\cdot
          \frac{\sigma_{30}\sigma_{31}}{\xi_0^2},&~
  f_{11}&=\frac{\xi_0^2}{\sigma_0^2}\cdot
          \frac{\sigma_{30}\sigma_{32}}{\xi_0^2}.
\endalignat $$
Then $H^0(P,2H_\sigma-2H_\xi)\simeq\Bbb C
      \langle f_1,\dots,f_{11}\rangle$
and  $\text{Bs}|2H_\sigma-2H_\xi|=\Bbb P_\xi^1\times\{(1:0:\dots:0)\}$.
For any $f=\sum_{i=1}^{11} c_if_i$ with $c_1\not=0$,
$\text{div}(f\cdot s_2)$ is nonsingular along
$\text{Bs}|2H_\sigma-2H_\xi|$,
thus any general member $M$ of $|2H_\sigma-2H_\xi|$ is smooth.
Putting $L=H_M$,we obtain a polarized manifold $(M,L)$ as desired.
In this case,$|L-H_\xi|$ makes $M$ a desingularization
of a variety of degree six in $\Bbb P^7$.

(3.24.5)
Even when $\Cal E\simeq\Cal O(1,1,3,3)$,
we can show the existence of $(M,L)$ similarly.
We have
$P\simeq\{(\xi_0:\xi_1)\times(\sigma_0:\sigma_1:
                              \sigma_{20}:\sigma_{21}:\sigma_{22}:
                              \sigma_{30}:\sigma_{31}:\sigma_{32})
  \in\Bbb P_\xi^1\times\Bbb P_\sigma^7
 |\xi_0:\xi_1=\sigma_{20}:\sigma_{21}
             =\sigma_{21}:\sigma_{22}
             =\sigma_{30}:\sigma_{31}
             =\sigma_{31}:\sigma_{32}\}$
and $H^0(P,2H_\sigma-2H_\xi)\simeq\Bbb C
      \langle f_1,\dots,f_{13}\rangle$
,where
  $$ \alignat 4
     f_1&=\frac{\xi_0^2}{\sigma_0^2}\cdot
          \frac{\sigma_0\sigma_{20}}{\xi_0^2},&~
     f_2&=\frac{\xi_0^2}{\sigma_0^2}\cdot
          \frac{\sigma_0\sigma_{30}}{\xi_0^2},&~
     f_3&=\frac{\xi_0^2}{\sigma_0^2}\cdot
          \frac{\sigma_1\sigma_{20}}{\xi_0^2},&~
     f_4&=\frac{\xi_0^2}{\sigma_0^2}\cdot
          \frac{\sigma_1\sigma_{30}}{\xi_0^2},\\
     f_5&=\frac{\xi_0^2}{\sigma_0^2}\cdot
          \frac{\sigma_{20}^2}{\xi_0^2},&~
     f_6&=\frac{\xi_0^2}{\sigma_0^2}\cdot
          \frac{\sigma_{20}\sigma_{21}}{\xi_0^2},&~
     f_7&=\frac{\xi_0^2}{\sigma_0^2}\cdot
          \frac{\sigma_{21}^2}{\xi_0^2},&~
     f_8&=\frac{\xi_0^2}{\sigma_0^2}\cdot
          \frac{\sigma_{20}\sigma_{30}}{\xi_0^2},\\
     f_9&=\frac{\xi_0^2}{\sigma_0^2}\cdot
          \frac{\sigma_{20}\sigma_{31}}{\xi_0^2},&~
  f_{10}&=\frac{\xi_0^2}{\sigma_0^2}\cdot
          \frac{\sigma_{21}\sigma_{31}}{\xi_0^2},&~
  f_{11}&=\frac{\xi_0^2}{\sigma_0^2}\cdot
          \frac{\sigma_{30}^2}{\xi_0^2},&~
  f_{12}&=\frac{\xi_0^2}{\sigma_0^2}\cdot
          \frac{\sigma_{30}\sigma_{31}}{\xi_0^2},\\
  f_{13}&=\frac{\xi_0^2}{\sigma_0^2}\cdot
          \frac{\sigma_{31}^2}{\xi_0^2}.
\endalignat $$
Since
 $\text{Bs}|2H_\sigma-2H_\xi|=\{\sigma_{20}=\sigma_{21}=\dots
                                    =\sigma_{32}=0~\text{in}~P\}$,
if we take $f=\sum_{i=1}^{13} c_if_i$ with $c_1c_4-c_2c_3\not=0$,
$\text{div}(f\cdot s_2)$ is nonsingular along
$\text{Bs}|2H_\sigma-2H_\xi|$.
Thus any general member $M$ of $|2H_\sigma-2H_\xi|$ is smooth.
Putting $L=H_M$,we obtain a polarized manifold $(M,L)$ as desired,
and $|L-H_\xi|$ makes $M$ a desingularization
of a variety of degree six in $\Bbb P^7$.

Summarizing the results in \S3,we obtain the following.

\proclaim{{\rm(3.25)}Theorem}
Let $(M,L)$ be a polarized manifold of the type (1.5.3).
Then $g(C)$,the genus of $C$,is 0 or 1,
$\Cal E:=\varPhi_*\Cal O_M[L]$ is a locally free sheaf on $C$,
$M\in|2H(\Cal E)+B_{\Bbb P(\Cal E)}|$
for some line bundle $B$ on $C$,and $L=[H(\Cal E)]_M$.
Putting $d=L^n,e=c_1(\Cal E)$,and $b=\deg B$,
we have the following results.

When $g(C)=1$,we have $1\leq d\leq 6,e=d-2,b=4-d$,and
\roster
\item"(i)" if $d=1$ or $2$,then $\Cal E$ is not ample;
\item"(ii)" if $d=3$ or $4$,then $\Cal E$ is ample as long as
            it is indecomposable;
\item"(iii)" if $d=5$ or $6$,then $\Cal E$ is ample.
\endroster

When $g(C)=0$,we have $C\simeq\Bbb P_\xi^1,1\leq d\leq 12,e=d-4,b=8-d,
                       M\in|2H(\Cal E)+bH_\xi|$,
and their lists are in the table below.

\newpage
$$\vbox{\offinterlineskip \eightpoint
  \halign{\strut\vrule# & \enskip\hfil#\hfil\enskip & \vrule#
                        & \enskip\hfil#\hfil\enskip & \vrule#
                        & \enskip\hfil#\hfil\enskip & \vrule# \cr
  \noalign{\hrule}
    & $d$ && $\Cal E$ && $(M,L)$  &\cr
  \noalign{\hrule}
    &  1  && $\Cal O(-3,0,0,0)$ && The existence is uncertain. &\cr
    &     && $\Cal O(-3,0,0,0,0)$ && $_{''}$  &\cr
    &     && $\Cal O(-2,-1,0,0)$ && $_{''}$  &\cr
    &     && $\Cal O(-2,-1,0,0,0)$ && $_{''}$  &\cr
    &     && $\Cal O(-1,-1,-1,0)$ && $_{''}$  &\cr
    &     && $\Cal O(-1,-1,-1,0,0)$ && $_{''}$  &\cr
  \noalign{\hrule}
    &  2  && $\Cal O(-2,0,0,0)$ && $_{''}$  &\cr
    &     && $\Cal O(-2,0,0,0,0)$ && $_{''}$  &\cr
    &     && $\Cal O(-1,-1,0,0)$ && $_{''}$  &\cr
    &     && $\Cal O(-1,-1,0,0,0)$ && $_{''}$  &\cr
  \noalign{\hrule}
    &  3  && $\Cal O(-2,0,0,1)$ && $_{''}$  &\cr
    &     && $\Cal O(-2,0,0,0,1)$ && $_{''}$  &\cr
    &     && $\Cal O(-1,-1,0,1)$ && $_{''}$  &\cr
    &     && $\Cal O(-1,-1,0,0,1)$ && $_{''}$  &\cr
    &     && $\Cal O(-1,0,0,0)$ && $_{''}$  &\cr
    &     && $\Cal O(-1,0,0,0,0)$ && $_{''}$  &\cr
  \noalign{\hrule}
    &  4  && $\Cal O(-1,0,0,1)$ && $_{''}$  &\cr
    &     && $\Cal O(-1,0,0,0,1)$ && $_{''}$  &\cr
    &     && $\Cal O(0,0,0,0)$ && $|L|$ makes $M$ a quadruple covering
                                  of $\Bbb P^3$.  &\cr
    &     && $\Cal O(0,0,0,0,0)$ && $|L|$ makes $M$ a quadruple covering
                                  of $\Bbb P^4$.  &\cr
  \noalign{\hrule}
    &  5  && $\Cal O(-1,0,0,2)$ && Bs$|L|$ is a point.  &\cr
    &     && $\Cal O(-1,0,1,1)$ && $_{''}$  &\cr
    &     && $\Cal O(-1,0,0,1,1)$ && $_{''}$  &\cr
    &     && $\Cal O(0,0,0,1)$ && $|L|$ makes $M$ the normalization of a
                                hypersurface of degree five in $\Bbb P^4$.&\cr
    &     && $\Cal O(0,0,0,0,1)$ && $|L|$ makes $M$ the normalization of a
                                 hypersurface of degree five in $\Bbb P^5$.&\cr
  \noalign{\hrule}
    &  6  && $\Cal O(-1,1,1,1)$ && $M$ is a double covering of
                                   $\Bbb P_\xi^1\times\Bbb P_\sigma^2$
                                   with branch locus being  &\cr
    &     &&                    && a smooth divisor of bidegree (4,2).
                                   $L=[H_\xi+H_\sigma]_M$. &\cr
    &     && $\Cal O(0,0,1,1)$ && Exist. &\cr
    &     && $\Cal O(0,0,0,1,1)$ && $_{''}$  &\cr
    &     && $\Cal O(0,0,0,2)$ && $_{''}$  &\cr
  \noalign{\hrule}
    &  7  && $\Cal O(0,0,1,2)$ && $_{''}$  &\cr
    &     && $\Cal O(0,1,1,1)$ && $_{''}$  &\cr
    &     && $\Cal O(0,0,1,1,1)$ && $_{''}$  &\cr
  \noalign{\hrule}
    &  8  && $\Cal O(0,1,1,1,1)$ && $M$ is a double covering of
                                   $\Bbb P_\xi^1\times\Bbb P_\sigma^3$
                                   with branch locus being  &\cr
    &     &&                    && a smooth divisor of bidegree (2,2).
                                   $L=[H_\xi+H_\sigma]_M$. &\cr
    &     && $\Cal O(0,1,1,2)$ && $M$ is a double covering of a divisor of
                                  bidegree (1,1) &\cr
    &     &&                   && on $\Bbb P_\xi^1\times\Bbb P_\sigma^3$.
                                  $L=[H_\xi+H_\sigma]_M$. &\cr
    &     && $\Cal O(1,1,1,1)$ && $M$ is a smooth divisor of bidegree (2,2)
                                  on $\Bbb P_\xi^1\times\Bbb P_\sigma^3$.
                                  $L=[H_\xi+H_\sigma]_M$. &\cr
  \noalign{\hrule}
    &  9  && $\Cal O(1,1,1,1,1)$ && $M$ is the blowing-up of
                                    $\Bbb P_\sigma^4$ with center
                                    being a complete  &\cr
    &     &&                     && intersection of two hyperquadrics.
                                    $L=[H_\xi+H_\sigma]_M$. &\cr
    &     && $\Cal O(1,1,1,2)$ && $M$ is the strict transform of a smooth
                                  hyperqubic in $\Bbb P_\sigma^4$ by the &\cr
    &     &&                   &&  blowing-up of $\Bbb P_\sigma^4$ with center
                                    being a $\Bbb P^2$.
                                    $L=[H_\xi+H_\sigma]_M$. &\cr
  \noalign{\hrule}
    & 10  && $\Cal O(1,1,1,1,1,1)$ && $M\simeq\Bbb P_\xi^1\times Q$,where
                                      $Q$ is a smooth hyperquadric in
                                      $\Bbb P_\sigma^5$.
                                      $L=[H_\xi+H_\sigma]_M$. &\cr
    &     && $\Cal O(1,1,1,1,2)$ && $M$ is the blowing-up of
                                    a hyperquadric in $\Bbb P_\sigma^5$
                                    with center &\cr
    &     &&                     && being a smooth quadric surface.
                                    $L=[H_\xi+H_\sigma]_M$. &\cr
  \noalign{\hrule}
}}$$
$$\vbox{\offinterlineskip \eightpoint
  \halign{\strut\vrule# & \enskip\hfil#\hfil\enskip & \vrule#
                        & \enskip\hfil#\hfil\enskip & \vrule#
                        & \enskip\hfil#\hfil\enskip & \vrule# \cr
  \noalign{\hrule}
    & $d$ && $\Cal E$ && $(M,L)$  &\cr
  \noalign{\hrule}
    & 10  && $\Cal O(1,1,2,2)$ && $M$ is a desingularization of a
                                  complete intersection of  &\cr
    &     &&                   && two hyperquadrics in $\Bbb P_\sigma^5$.
                                  $L=[H_\xi+H_\sigma]_M$. &\cr
    &     && $\Cal O(1,1,1,3)$ && $_{''}$  &\cr
  \noalign{\hrule}
    & 11  && $\Cal O(1,2,2,2)$ && $|L-H_\xi|$ makes $M$ a desingularization
                                  of a three-  &\cr
    &     &&                   && dimensional variety of degree five
                                  in $\Bbb P^6$. &\cr
  \noalign{\hrule}
    & 12  && $\Cal O(1,1,3,3)$ && $|L-H_\xi|$ makes $M$ a desingularization
                                  of a three-  &\cr
    &     &&                   && dimensional variety of degree six
                                  in $\Bbb P^7$. &\cr
    &     && $\Cal O(1,2,2,3)$ && $_{''}$  &\cr
    &     && $\Cal O(2,2,2,2)$ && $M\simeq\Bbb P_\xi^1\times
                                          \Bbb P_\mu^1\times
                                          \Bbb P_\lambda^1$ and
                                  $L=2H_\xi+H_\mu+H_\lambda$. &\cr
  \noalign{\hrule}
}}$$
\endproclaim

\subhead
\S4 The case of a Veronese fibration over a curve
\endsubhead

In this section we study the case (2) of the theorem (1.5),
using the argument in \cite{F;(II.13.10)}.

(4.1)
Put $H=K+2L$,then $\Cal E:=\varPhi_*\Cal O_M[H]$
is a locally free sheaf of rank three on $C$ and
$(M,H)$ is the scroll of $\Cal E$.
We have $L=2H+\varPhi^*B$ for some $B\in\text{Pic}(C)$.
Similarly as before,we put $d=L^3,e=c_1(\Cal E),b=\deg B$
and denote by $g(C)$ the genus of $C$.
Then $e\geq 0,e+b=1$,and $d=8e+12b$.
By the canonical bundle formula,we obtain that
$K^C+\det\Cal E+2B=0$,hence $2g(C)-2+e+2b=0$.
{}From these results,$(e,d)=(0,12)$ or $(2,4)$.

(4.2)
When $(e,d)=(0,12)$,we have $b=1$ and $g(C)=0$,
hence $C\simeq\Bbb P^1,B=\Cal O(1)$,and
$\Cal E\simeq\Cal O(e_1)\oplus\Cal O(e_2)\oplus\Cal O(e_3)$
for $e_1,e_2,e_3\in\Bbb Z$.
For each $1\leq i\leq 3$,a natural surjection $\Cal E\to\Cal O(e_i)$
gives a section $Z_i$ of $\varPhi$ and $H_{Z_i}=\Cal O(e_i)$.
Since $e_1+e_2+e_3=e=0$ and $L_{Z_i}=\Cal O(2e_i+1)$ is ample,
we have $e_1=e_2=e_3=0$ and $\Cal E\simeq\Cal O_C^{\oplus3}$,
thus $M\simeq\Bbb P_\xi^1\times\Bbb P_\sigma^2$
and $L=H_\xi+2H_\sigma$.
We note that this $(M,L)$ is already obtained in (2.9.3).

(4.3)
When $(e,d)=(2,4)$,we have $b=-1$ and $g(C)=1$.
Hence $C$ is an elliptic curve and $\det\Cal E+2B=0$
since $K^C=\Cal O_C$.
Let $Q$ be any quotient bundle of $\Cal E$.
If rank~$Q=1$,then $Z:=\Bbb P_C(Q)$ is a section of
$\varPhi$ and $HZ=c_1(Q)$.
Then $c_1(Q)\geq 1$ since $0<LZ=2c_1(Q)-1$.
If rank~$Q=2$,then $D:=\Bbb P_C(Q)\in|H-\varPhi^*\Cal F|$,
where $\Cal F$ is the kernel of $\Cal E\to Q$.
Since $0<L^2D=4(1-c_1(\Cal F))$,we have $c_1(Q)=e-c_1(\Cal F)\geq 2$.
In both cases we have $(\text{rank}~Q)\cdot c_1(\Cal E)
                      <(\text{rank}~\Cal E)\cdot c_1(Q)$,
hence $\Cal E$ is stable.
Conversely,let $\Cal E$ be a semistable vector bundle on $C$
with rank~$\Cal E=3$ and $c_1(\Cal E)=2$.
We put $M=\Bbb P_C(\Cal E),H=H(\Cal E)$ and let
$\varPhi:M\to C$ be the bundle map.
By the semistability criterion in \cite{Mi;(3.1)},
$3H-\varPhi^*(\det\Cal E)$ is nef.
Since $C$ is an elliptic curve,we can find some $B\in\text{Pic}(C)$
satisfying $\det\Cal E+2B=0$.
Then $3(2H+\varPhi^*B)=2(3H+\varPhi^*(2B))-\varPhi^*B$ is ample.
Hence $L:=2H+\varPhi^*B$ is ample and
$(M,L)$ is a polarized manifold of the type (1.5.2).

(4.4)
Summing up,we obtain the following theorem.
\proclaim{Theorem}
Let $(M,L)$ be a polarized manifold of the type (1.5.2).
We put $d=L^3$ and denote by $g(C)$ the genus of $C$.
Then $(M,L)$ is one of the following two types.
\roster
\item"({\rm I})"  $g(C)=0$,hence $C\simeq\Bbb P_\xi^1$;
                 $d=12,M\simeq\Bbb P_\xi^1\times\Bbb P_\sigma^2$,
                 and $L=H_\xi+2H_\sigma$.
\item"({\rm II})" $g(C)=1$ and $M\simeq\Bbb P_C(\Cal E)$,
                 where $\Cal E:=\varPhi_*\Cal O_M[K+2L]$
                 is a stable vector bundle of rank three on $C$
                 with $c_1(\Cal E)=2$;
                 $d=4$ and $L=2H(\Cal E)+\varPhi^*B$,
                 where $B\in\text{Pic}(C)$ with $\det\Cal E+2B=0$.
\endroster
\endproclaim

\subhead
\S5  Remaining cases
\endsubhead

In this section we study the cases (1) and (5) of the theorem (1.5).

(5.1)
We first study the case (1.5.5).
This case is a kind of ``general type'' and we mainly study
$(M,L)$ of small $\varDelta$-genus.
Since $L$ is ample,$0\leq(K+(n-2)L)L^{n-1}=2g-2-L^n=4-d$,
where $d=L^n$.
Hence we have $1\leq d\leq 4$.
For any polarized manifold,$\varDelta=0$ implies $g=0$
(\cite{F1;(1.9)} and \cite{F2;(4.1)}).
Thus here we have $\varDelta\geq 1$.
The case $d=1$ is difficult and yet to be studied.

(5.2)
Suppose that $d=2$.
When $\varDelta=1$,$|L|$ makes $M$ a double covering of
$\Bbb P^n$ with branch locus being a smooth hypersurface
of degree eight (\cite{F3;(2.5)}).
When $\varDelta=2$,Bs$|L|$ is a finite set
(\cite{F1;(1.9)} and \cite{F5;(1.17)}).
When $\varDelta\geq 3$,we do not have any satisfactory result.

(5.3)
Suppose that $d=3$.
Then we have $\varDelta\not=1$ by \cite{F1;(1.9)} and \cite{F2;(4.1)}.
When $\varDelta=2$,dim Bs$|L|\leq 0$
(\cite{F1;(1.9)} and \cite{F5;(1.14.5)}).
If $\varDelta=2$ and Bs$|L|=\phi$,$|L|$ makes $M$
a triple covering of $\Bbb P^n$.
If $\varDelta=2$ and Bs$|L|\not=\phi$,
we have the following results in \cite{F9}.
\roster
\item"a)" Bs$|L|$ consists of one simple point $p$.
\item"b)" Let $\pi:M_1\to M$ be the blowing-up at $p$ and
          let $E$ be the exceptional divisor over $p$.
          Then Bs$|\pi^*L-E|=\phi$.
\item"c)" $|\pi^*L-E|$ gives a morphism $\rho:M_1\to\Bbb P^n$
          of degree two.Every fiber $X$ of $\rho$ with $\dim X>0$
          is an irreducible curve such that $EX=1$.
\item"d)" $\rho(E)$ is a hyperplane in $\Bbb P^n$.
\endroster
When $\varDelta\geq 3$,we do not have any satisfactory result.

(5.4)
Suppose that $d=4$.
Then we have $(K+(n-2)L)L^{n-1}=0$.
Since $K+(n-2)L$ is nef,by the base point free theorem
(cf. \cite{KMM}) and \cite{F5;Appendix},there is a fibration
$f:M\to V$ and an ample line bundle $A$ on $V$ such that
$K+(n-2)L=f^*A$.
Thus we have $K+(n-2)L=0$ since $(f^*A)L^{n-1}=0$.
Then $\varDelta(M,L)=2$ by \cite{F4;(1.11)},
hence dim Bs$|L|\leq 1$ by \cite{F1;(1.9)}.
When dim Bs$|L|=1$,$n\leq 4$ by \cite{F5;(1.17)},
and we have the following results by \cite{F5;(1.14)\&(2.4)}.
\roster
\item"a)" $Y:=\text{Bs}|L|$ is a smooth rational curve.
\item"b)" Let $\pi:M'\to M$ be the blowing-up of $Y$ and
          let $E$ be the exceptional divisor over $Y$.
          Then Bs$|\pi^*L-E|=\phi$.
\item"c)" Let $W$ be the image of the morphism $M'\to\Bbb P^{n+1}$
          defined by $|\pi^*L-E|$.Then $\dim W=n-1,\deg W=3$
          and $\varDelta(W,\Cal O_W(1))=0$.
\item"d)" $E$ is a section of the morphism $\rho:M'\to W$.
\item"e)" $\rho$ is flat and every fiber of $\rho$ is an
          irreducible curve of arithmetic genus one.
\item"f)" $M'$ is a double covering of a $\Bbb P^1$-bundle
          $V:=\Bbb P(\Cal O_E\oplus[-2E]_E)$ over $E\simeq W$.
          The image $S$ of $E$ by the morphism $\beta:M'\to V$
          is the unique member of $|H_\zeta-[-2E]_V|$,
          where $H_\zeta$ is the tautological bundle on $V$.
          The branch locus $B$ of $\beta$ is $B^*+S$,
          where $B^*$ is a smooth member of $|3H_\zeta|$
          with $B^*\cap S=\phi$.
\endroster
When dim Bs$|L|\leq 0$,we have Bs$|L|=\phi$ by \cite{F2;(4.1)}.
Let $\rho$ be the morphism $M\to\Bbb P^{n+1}$ defined by $|L|$.
We put $W=\rho(M)$ and $w=\deg W$,then $4=L^n=w\cdot\deg\rho$.
Hence $(\deg\rho,w)=(1,4)$ or (2,2).
In the former case,$\rho$ is birational and moreover
$M\simeq W$ by \cite{F;(10.8.1)}.
In the latter case,$\rho$ is a double covering of a hyperquadric $W$.
Furthermore $W$ turns out to be smooth,and the branch locus of $\rho$
is a smooth hypersurface section and is connected (\cite{F;(10.8.2)}).

(5.5)
Finally we study the case (1) of the theorem (1.5).
We use the theory of minimal reduction in \cite{F7;(1.9)}
and \cite{F;(11.11)}.
Clearly $M$ is the blowing-up of another manifold $M_1$ at one point
and $E$ is the exceptional divisor with $LE=1$.
Moreover there is an ample line bundle $L_1$ on $M_1$ such that
$L+E$ is the pullback of $L_1$.
In such a case we say that $(M,L)$ is a simple blow-up of $(M_1,L_1)$.
By the above theory,we obtain a sequence of simple blow-ups
  $$ (M,L)\to(M_1,L_1)\to\dots\to(M_r,L_r)=(M',L') $$
with the following properties:
\roster
\item"a)" $K+(n-1)L=[K'+(n-1)L']_M$ for the canonical bundle $K'$ of $M'$,
          hence $K'+(n-1)L'$ is nef and $g(M',L')=g(M,L)=3$;
\item"b)" $L^n=(L')^n-r$;
\item"c)" $(M',L')$ is not of the type (1.5.1) i.e. $(M',L')$ is minimal.
\endroster
If $(M',L')$ is of the type (1.5.3) or (1.5.4),
then we can derive a contradiction as in \cite{F7;(1.9)}.
Thus $(M',L')$ is of the type (1.5.2) or (1.5.5).

(5.6)
When $(M',L')$ is of the type (1.5.2),
we can apply the argument in \S4.
If $M'\simeq\Bbb P_\xi^1\times\Bbb P_\sigma^2$ and $L'=H_\xi+2H_\sigma$,
then we can find a curve $Z$ on $M'$ with $L'Z=1$.
This is a contradiction by \cite{F7;(1.9)},
thus $(M',L')$ is of the type (4.4.II).
Since $4=(L')^3=L^3+r$,the number of points at which $M'$ is blown up
is less than four.

(5.7)
When $(M',L')$ is of the type (1.5.5),
$0\leq(K'+(n-2)L')(L')^{n-1}=4-(L')^n$ since $K'+(n-2)L'$ is nef.
Then we have $2\leq L^n+r=(L')^n\leq 4$,hence
$(L^n,r,(L')^n)=(1,1,2),(1,2,3),(1,3,4),(2,2,4)$,or (3,1,4).

\enddocument